\documentclass{article}
\usepackage{amsmath,amsfonts,amssymb,amsthm,graphicx}
\usepackage{dsfont}    
\usepackage{enumerate}
\usepackage{color}
\usepackage[numbers]{natbib}
\usepackage{multirow}
\usepackage{comment}
\usepackage{booktabs,subcaption,dcolumn}

\usepackage{algorithm}
\usepackage[noend]{algpseudocode}  
\usepackage[colorlinks=true,citecolor=blue,pdfpagemode=UseNone,pdfstartview=FitH]{hyperref}

\newif\ifFULL
\ifFULL

\fi

\newif\ifRuodu
\ifRuodu
  \renewcommand{\ge}{\geqslant}
  \renewcommand{\le}{\leqslant}
  \renewcommand{\epsilon}{\varepsilon}
  
\fi

\renewcommand{\d}{\,\mathrm{d}}

\newcommand{\p}{\mathbb{P}}

\newcommand{\E}{\mathbb{E}}    
\newcommand{\var}{\mathrm{var}}  
\newcommand{\R}{\mathbb{R}}    
\newcommand{\N}{\mathbb{N}}    

  \newcommand{\id}{\mathds{1}}

\usepackage{setspace}

\theoremstyle{plain}
\newtheorem{theorem}{Theorem} 

\newtheorem{lemma}{Lemma}
\newtheorem{proposition}{Proposition}

\theoremstyle{definition}

\newtheorem{example}{Example}

\theoremstyle{remark}
\newtheorem{remark}{Remark}

\renewcommand{\cite}{\citet}  
\DeclareMathOperator*{\argmax}{arg\,max}

\title{Testing   mean and variance  by e-processes} 
\author{Yixuan Fan\thanks{
Center for Applied Mathematics,
  Tianjin University, 
  Tianjin, China
  E-mail: \href{mailto:yxfanjlu@gmail.com}{yxfanjlu@gmail.com}}
\and Zhanyi Jiao%
\thanks{Department of Mathematics,
Illinois State University,
  Normal, Illinois, USA.
  E-mail: \href{mailto:zjiao1@ilstu.edu}{zjiao1@ilstu.edu}}  \and
  Ruodu Wang\thanks%
  {Department of Statistics and Actuarial Science,
  University of Waterloo,
  Waterloo, Ontario, Canada.
  E-mail: \href{mailto:wang@uwaterloo.ca}{wang@uwaterloo.ca}.}}

\begin{document}
\maketitle
 
\begin{abstract}
We address the problem of testing conditional mean and conditional variance for non-stationary data. We build e-values and p-values for four types of non-parametric composite hypotheses with specified mean and variance as well as other conditions on the shape of the data-generating distribution. These shape conditions include symmetry, unimodality, and their combination. Using the obtained e-values and p-values, we construct tests via e-processes, also known as testing by betting, as well as some tests based on combining p-values for comparison. Although we mainly focus on one-sided tests, the two-sided test for the mean is also studied. Simulation and empirical studies are conducted under a few settings, and they illustrate features of the methods based on e-processes. \vspace{2ex}
\end{abstract}

\section{Introduction}

Testing mean and variance in various settings is a classic problem in statistics. In parametric inference concerning testing the mean, well-known tests like Student's t-test and z-test, as well as tests related to variance such as the chi-squared test and the F-test, are commonly employed; see e.g., \cite{LRC86}. Parametric tests always come with   assumptions about the forms of the population distribution from which samples are derived. Deviating from these assumptions can lead to significantly flawed results. For situations where these assumptions might be compromised, non-parametric methods provide a great alternative. 
Certainly, non-parametric methods may also make strong assumptions on the underlying population, such as finite or bounded moments, but  not on the specific parametric forms.
Comprehensive and well-established methods of non-parametric techniques for testing means and variances can be found in e.g., \cite{C99} and \cite{HWC13}. 
Different from the classic settings, we will consider the problem of testing composite hypotheses in which data are not stationary. 

Suppose that a tester has sequentially arriving, possibly dependent, data points $X_1,X_2,\dots$, each from an unknown distribution, possibly different. The tester is interested in testing whether
\begin{align}\label{eq:intro}
\E[X_i|\mathcal F_{i-1}]\le \mu_i \mbox{~~and~~} 
\var(X_i|\mathcal F_{i-1})\le \sigma_i^2 \mbox{~~for each $i$},
\end{align}
where $\mathcal F_{i-1}$  is the $\sigma$-algebra generated by $X_1,\dots,X_{i-1}$, and $\mu_i$ and $\sigma_i$  are $\mathcal F_{i-1}$-measurable. 
All conditional expectations are in the almost sure sense. 
If independence is further assumed, then this problem reduces to the classic problem of testing mean and variance.
Testing conditional mean and conditional variance is common in some contexts such as forecasting (e.g., \cite{HZ22}) and financial risk assessment (e.g., \cite{FZ16}).

Problem \eqref{eq:intro} can be interpreted in two different ways, omitting ``conditional" here:
\begin{enumerate}
\item[(A)] testing both the   mean and the   variance;
\item[(B)] testing the mean under the knowledge of an upper bound on the   variance. 
\end{enumerate}
The interpretation (A) is relevant when the tester is interested in whether a time-series has switched away from a given regime with specified mean and variance bounds. 
We mainly use interpretation (A), while keeping in mind that  interpretation (B) is   useful when comparing with the literature. Of course, one could also interpret  \eqref{eq:intro}  as testing the variance under the knowledge of an upper bound on the mean.

Clearly,   problem \eqref{eq:intro} is a composition of many complicated, non-parametric, composite hypotheses on each observation.
The key challenge in this setting is that the data points are not iid, and hence we cannot  make inference of the distributions themselves.

This problem can be addressed with the following general methodology, called e-testing or testing by betting, a successful example being \cite{WR24}.
We first consider a simpler problem: constructing an e-value from  one random variable from each data point with the corresponding hypothesis on its mean and variance, which corresponds to $n=1$. 
 For a general background  on e-values in hypothesis testing; see  \cite{VW21}, \cite{GDK20}, and the review by \cite{RGVS22}.
After obtaining these e-values, 
we combine them, usually by forming an e-process,  to construct a test for the overall   hypothesis. 
Alternatively, we can construct p-values instead of e-values, but the power of such a strategy is usually quite weak, as seen from our experiments. 

Section \ref{sec:2} formally describes the hypotheses and defines e-variables, e-processes, and p-variables. 
As mentioned above, we will first  address the case of one data point, i.e., $n=1$, presented in Section \ref{sec:single}. 
We consider four types of composite hypotheses on  mean, variance and the shape of the distribution: symmetry, unimodality and their combination. 
Our main results are ways that are  \emph{optimal}, in a natural sense, to constructions of p-values and e-values in this setting.   
Although our main methodology is based on e-processes, we present results also for p-values, which may be useful in multiple testing, not treated in this paper; for instance, p-values are the inputs of the standard procedure of \cite{BH95}. Considering a non-parametric composite hypothesis with a given mean and variance as the baseline case, assuming symmetry approximately improves the baseline p-variable by a multiplicative factor of 1/2, unimodality by a factor of 4/9, and both by a factor of 2/9. Similarly, the corresponding baseline e-variable is improved by multiplicative factors of 2, 1, and 2, respectively, in these scenarios; recall that smaller p-values are more useful, whereas larger e-values are more useful.

We propose in Section \ref{sec:multi} several methods to test using multiple data points, thus addressing the main task of the tester. The main proposals are e-process based tests, which follow the idea of testing by betting in \cite{S21}, \cite{WRB20}  and \cite{WR24}.  Although we mainly focus on  one-sided hypotheses, 
our methodology can  be easily adapted to test  the two-sided hypothesis on the mean,   that is, 
 $$
\E[X_i|\mathcal F_{i-1}]\in [ {\mu}^L_i, {\mu}^U_i]  \mbox{~~and~~} \var(X_i|\mathcal F_{i-1})\le \sigma_i^2 \mbox{~~for each $i$},
 $$
where $[ {\mu}^L_i, {\mu}^U_i]$  is an interval or a singleton  for each $i$; this   is discussed in Section \ref{sec:two-side}.

The closest methodological work related to this paper is \cite{WR24}, where the authors test in a non-parametric setting the conditional mean of sequential data, which are assumed to be bounded within a pre-specified range, thus a generally smaller class of distributions.
Our problem and methodology are different from \cite{WR24} in the sense that we assume a bounded variance instead of a bounded range. Since a bounded range implies bounded variance, the assumption needed to apply our methodology is weaker than in the setting of \cite{WR24}, following interpretation (B) of the main testing problem. 
 Moreover, we are able to utilize the additional information on the distributional shape to obtain better e-values than without such information. 
 A great advantage of the tests of \cite{WR24} is that their power adapts to the unknown true variance of the distribution if data come from an iid population. Our method based on the growth rate of empirical e-values has a similar feature, which uses a   betting strategy similar to that of \cite{WR24}. Another closely related methodology is \cite{WWZ22}, where, other statistical functions are tested other than the mean. Once e-variables are constructed, we will build e-processes in a similar way to \cite{WWZ22}.   The methods  of  \cite{HRMS20, HRMS21} and \cite{WR23}   based on exponential test supermartingales---exponential processes that form supermartingales with initial value one---which are e-processes, 
 can also be applied to test \eqref{eq:intro}. These methods differ from ours as our e-process is obtained by combining individual e-variables.

Section \ref{sec:simulation} provides simulation studies for the proposed methods and compare them with the method of \cite{WR24} when the model has both bounded support and bounded variance
and with methods based on exponential test supermartingale of \cite{HRMS21} and \cite{WR23}.
Section \ref{sec:emp} contains empirical studies using financial asset return data during the 2007--2008 financial crisis, further demonstrating the effectiveness of the e-process based methods. 
Section \ref{sec:7} concludes the paper. All proofs in the paper are provided in the Supplementary Material.

\section{General setting}
\label{sec:2}

\subsection{Hypotheses to test}

We first describe our main testing problem. 
Let $n$ be a positive integer or $\infty$, and denote by $[n]=\{1,\dots,n\}$. 
Throughout, fix a sample space.
Suppose that  data points  $(X_i)_{i \in [n]}$ arrive sequentially, each possibly from a different distribution,
and not necessarily independent. 
  A hypothesis is a collection $H$ of probability measures that govern $(X_i)_{i \in [n]}$.
Denote by   $\mathcal F_{i}$ the $\sigma$-field generated by $X_1,\dots,X_i$ for $i\in[n]$ with $\mathcal F_0$ being the trivial $\sigma$-field. 
The main hypotheses of interest  are variations (by adding shape information) of the following hypothesis 
\begin{align}
H =\left \{Q: \E^Q[X_i|\mathcal F_{i-1}]\le \mu_i \mbox{~and~}\var^Q(X_i|\mathcal F_{i-1}) \le \sigma_i^2 \mbox{ for  $i\in [n]$}\right\},\label{eq:maintest}
\end{align}
where $\mu_i$ and $\sigma_i$ are $ \mathcal F_{i-1}$-measurable  for each $i\in [n]$; that is, they can be data-dependent on past observations. 
A simple case is 
\begin{align}
H =\left \{Q: \E^Q[X_i|\mathcal F_{i-1}]\le \mu \mbox{~and~}\var^Q(X_i|\mathcal F_{i-1}) \le \sigma^2 \mbox{ for  $i\in [n]$}\right\},\label{eq:maintest2}
\end{align}
where $\mu$ and $\sigma$ are two constants; that is, we would like to test whether data exhibit conditional mean and conditional variance in $(-\infty,\mu]\times [0,\sigma^2]$. Although \eqref{eq:maintest2} looks simpler, it is indeed equivalent to \eqref{eq:maintest} by noting that $\mu_i$ and $\sigma_i$ are $\mathcal F_{i-1}$-measurable, they can be absorbed into $X_i$ by considering $(X_i-\mu_i)/\sigma_i$ instead of $X_i$. Therefore, we will focus on the formulation \eqref{eq:maintest2} for the rest of the paper.  
If data are independent, but not necessarily identically distributed, then the problem is to test the unconditional mean and variance. 
We sometimes omit $Q$ in $\E^Q$ and $\var^Q$ when it is clear.

We will further consider hypotheses with additional shape information, by assuming that some, or all of the distributions of $X_1,\dots,X_n$ are unimodal, symmetric, or both. 
Below, all terms like ``increasing" and ``decreasing" are in the non-strict sense. 
A distribution on $\R$ is \emph{unimodal} if there exists $x\in \R$ such that 
the distribution has an   increasing density on $(-\infty, x)$ and a decreasing density on $(x,\infty)$;
it may have a point-mass at $x$. 
A distribution on $\R$ with mean $\mu$ is \emph{symmetric} if  for all $x\in \R$ it assigns equal probabilities to $(-\infty, \mu-x]$ and $[\mu+x,\infty)$. 
If a distribution with mean $\mu$ is both unimodal and symmetric, then its mode must be either  $\mu$ or an interval centered at $\mu$.

\begin{remark} \label{rem:R1-WR}
\label{rem:variance} The main question in \cite{WR24} is to test the  conditional mean $m$ with data taking values in $[0,1]$.
We note that any random variable with mean at most $m$ and range $[0,1]$ has variance at most $1/4$ (if $m\ge 1/2$) or $m(1-m)$ (if $m<1/2$), attained by a Bernoulli random variable. Therefore, our hypothesis with $\mu=m$ and $\sigma^2=1/4$ or $\sigma^2=m(1-m)$ has less restrictive assumptions than their setting (except they formulated two-sided hypotheses; see Remark \ref{rem:WR} below) and in particular, our setting can handle unbounded data.
\end{remark}

 \begin{remark}\label{rem:WR}
Our hypotheses are formulated as one-sided on both $\mu$ and $\sigma^2$. 
Certainly, all validity results remain true for the two-sided hypotheses.
Testing $\E^Q[X_i]\ge \mu$ is symmetric to testing $\E^Q[X_i]\le \mu$, but such symmetry does not hold for testing the variance. 
Building e-processes to test the two-sided hypothesis on the mean is discussed in Section \ref{sec:two-side}.
\end{remark} 

\subsection{P-variables and e-variables}
We formally define p-variables and e-variables, following 
   \cite{VW21}.  
A \emph{p-variable} $P$ for a hypothesis $H$ is a  random variable  that satisfies $Q(P\le \alpha)\le \alpha$ for all $\alpha \in (0,1)$ and all $Q\in H$. 
In other words, a p-variable is stochastically larger than $\mathrm {U}[0,1]$, often truncated at $1$. An \emph{e-variable} $E$ for a hypothesis $H$ is a $[0,\infty]$-valued random variable satisfying $\E^Q[E]\le1$ for all $Q\in H$.
E-variables are often obtained from stopping an \emph{e-process} $(E_t)_{t \geq 0}$, which is a non-negative stochastic process adapted to a pre-specified filtration, $(\mathcal F_i)_{i\in [n]}$ in our problem, such that $\mathbb{E}^Q[E_\tau] \leq 1$ for any stopping time $\tau$ and any $Q\in H$.

Some p-variables and e-variables are useless, like $P=1$ or $E=1$. 
A p-variable $P$ for $H$ is \emph{precise}  if  $\sup_{Q\in H}Q(P\le \alpha)=\alpha$ for each $\alpha\in (0,1)$, and an e-variable $E$ for $H$ is \emph{precise}  if $\sup_{Q\in H}\E^Q[E] = 1$. In other words, a p-variable or an e-variable being precise means that it is not wasteful in a natural sense. 
For instance, if  $\sup_{Q\in H}\E^Q[E] <1$, then we can multiply $E$ by a constant  larger than $1$.
Some imprecise e-variables may also be useful, such as those built on the Hoeffding inequality; see 
\cite{H63},  \cite{HRMS21} and \cite{WR24}.

A p-variable $P$ is \emph{semi-precise} for $H$ if  $\sup_{Q\in H}Q(P\le \alpha)=\alpha$ for each $\alpha\in (0,1/2]$. Semi-precise p-variables require the sharp probability bound $\sup_{Q\in H}Q(P\le \alpha)=\alpha$ only for  the  case $\alpha \le 1/2$  which is relevant for testing purposes. 
We will see that for some hypotheses, precise p-variables do not exist unless we rely on external randomization, but semi-precise ones do exist. 

Realizations of p-variables and e-variables are referred to as p-values and e-values. As is customary in the literature, we sometimes, but never in mathematical statements, use the two terms ``e-value" and ``e-variable" interchangeably. 
 
\section{Best p-  and e-variables for one data point}\label{sec:single}

\subsection{Setting}
We begin by considering the simple setting where one data point $X$ is available,
from which we build a p-variable or e-variable for the hypothesis. 
Although it may be unconventional to test based on one observation, there are several situations where this  construction becomes useful.
\begin{enumerate}
\item Testing by betting:  To construct an e-process, one needs to sequentially obtain  one e-value from each observation, or a batch of observations. This is the main setting in the current paper.
\item Testing multiple hypotheses: One observation is obtained for each  hypothesis, and  p-values or e-values for each of them are computed and fed into a multiple testing procedure such as that of \cite{BH95};  this setting is particularly relevant for the procedure of \cite{WR22} based on e-values, which yields false discovery rate control under arbitrary dependence.
Even if for some hypotheses there is only one data point, a p-value or e-value, even moderate, say $e=0.8$ or $e=1.2$, from this hypothesis may be useful for the overall testing problem; see \cite{IWR24} where e-values are used as weights, so $e=0.8$ or $e=1.2$ matters. 
\item Testing a global null: One may   first obtain a p-value or e-value for each experiment and then combine them to test the global null, as in meta-analysis; see \cite{VW20, VW21} and the references therein. 
\end{enumerate}
E-values are relevant for all of the three contexts, and p-values are relevant for the second and the third contexts.

We will focus on p-variables, which are decreasing functions of $X$, and e-variables, which are increasing functions of $X$.
Thus, a larger value of $X$ indicates stronger evidence against the null; this is intuitive because we are testing the mean less or equal to $\mu$ in \eqref{eq:maintest2}. 
This assumption on p-variables and e-variables will be made throughout the rest of the paper.

\begin{remark} In the contexts of multiple testing and sequential e-values, the dependence among several e-values or p-values obtained is preserved from the dependence among the data points, 
if the monotonicity assumption above holds. This will be helpful  when applying statistical methods based on dependence assumptions; see \cite{BY01} for the BH (\cite{BH95}) procedure with positive dependence and 
\cite{CRW22} for BH with negative dependence. Both concepts of dependence are preserved under monotone transforms. 
\end{remark}

\subsection{Two technical lemmas}

The following lemma establishes that the infimum of p-variables based on the same data point $X$ is still a p-variable. 
This result relies on  our assumption that p-variables are decreasing functions of $X$.
\begin{lemma}\label{lem:1} 
For a given observation $X$ and  hypothesis $H$,
the infimum of p-variables, which are assumed to be decreasing functions of $X$,  is a p-variable. As a consequence,
 there exists a smallest  p-variable.
\end{lemma}

Although the smallest p-variable  for $H$ exists, it may not be precise. Indeed, in Theorems \ref{th:2} and \ref{th:4} below we will see that there may not exist any precise p-variable for some hypotheses.
 
The following lemma allows us to convert conditions on distribution functions into conditions on the corresponding quantile functions.
For a probability measure $Q $, 
denote by 
$$T^Q_Y(\alpha) =\inf \{x\in \R: Q(Y\le x)\ge \alpha\} \mbox{~~~for $\alpha \in (0,1)$};$$
that is, $T^Q_Y$ is the left-quantile function of $Y$ under $Q$.
 \begin{lemma}\label{lem:2}
For a random variable $P$ and a hypothesis $H$,
\begin{enumerate}[(i)]
\item $P$  is   a p-variable  if and only if
$
\inf_{Q\in H}  T^Q_P (\alpha)\ge \alpha 
$
for all $\alpha \in(0,1)$;
\item $P$ is   a precise p-variable if and only if $
\inf_{Q\in H} T^Q_P (\alpha) = \alpha 
$
for all $\alpha \in(0,1)$;
\item $P$  is  a semi-precise p-variable if and only if $
\inf_{Q\in H} T^Q_P (\alpha) = \alpha 
$
for all $\alpha \in(0,1/2)$ and $
 \inf_{Q\in H}  T^Q_P (\alpha)\ge \alpha 
$ for $\alpha \in [1/2,1)$.
\end{enumerate}
\end{lemma}
The proof of Lemma \ref{lem:2} is essentially identical to that of Lemma 1 of \cite{VW20}, which gives the equivalence between probability statements and quantile statements for merging functions of p-values.
Our construction for precise and semi-precise p-variables will be based on computing $\alpha \mapsto \sup_{Q\in H} T_X^Q(1-\alpha)$ and its inverse function.

\subsection{Main results}

Recall that we have only one observation, denoted by $X$. 
We consider the following four classes of non-parametric composite hypotheses, where 
$\mu \in \R$ and $\sigma >0$.
 \begin{align*}
H(\mu,\sigma)&=\left \{Q:\E^Q[X]\le \mu \mbox{~and~}\var^Q(X) \le \sigma^2 \right\};
 \\ 
H_{\rm S}(\mu,\sigma)&=\{Q\in H(\mu,\sigma): \mbox{$X$ is symmetrically distributed}\};
 \\ 
H_{\rm U}(\mu,\sigma)&=\{Q\in H(\mu,\sigma): \mbox{$X$ is unimodally distributed}\};
 \\
H_{\rm US}(\mu,\sigma)&=H_{\rm U} (\mu,\sigma)\cap H_{\rm S}(\mu,\sigma).
 \end{align*}

For our main results on the ``best"  p-variables and e-variables, it will be clear from our proofs that the condition $\var^Q(X) \le \sigma^2$ in each hypothesis can be replaced by $\var^Q(X) = \sigma^2$,
and the condition $\E^Q[X] \le  \mu$ in each hypothesis can be replaced by $\E^Q[X] = \mu$.
All results remain true with any combinations of the above alternatives.
Possible improvement for the two-sided test is discussed in Section \ref{sec:two-side}.

The above four sets of distributions are studied in a very different context by \cite{LSWY18} to compute worst-case risk measures under model uncertainty  in finance. Some of our techniques for constructing p-variables use results from  \cite{LSWY18} and \cite{BKV20} for finding bounds on quantile, which is called the Value-at-Risk in finance.
 
In what follows, for $x\in \R$, we write $x_+=\max\{x,0\}$, 
$x_-=\max\{-x,0\}$, $x_+^2 =(x_+)^2 $, and $x_-^2=(x_-)^2$. We first consider the simplest case of testing $H(\mu,\sigma)$.
\begin{theorem}\label{th:1}
A precise p-variable for $H(\mu,\sigma)$ is  $P=(1+(X-\mu)_+^2/\sigma^2)^{-1}$,
and a precise   e-variable for $H(\mu,\sigma)$ is  $E=(X-\mu)_+^2/\sigma^2$.
\end{theorem}

Theorem \ref{th:1} can be seen as consequence of Cantelli's inequality.
It may be interesting to compare $P$ and $1/E$ obtained from Theorem \ref{th:1}.
Note that any e-variable can be converted into a p-variable  via the so-called calibrator $e\mapsto \min\{1/e,1\}$; see e.g., \cite{VW21}; this is an immediate consequence of Markov's inequality. 
As $1/E$ is a p-variable for an e-variable $E$, we have $P\le1/E$. 
In Theorem \ref{th:1}, we obtain $ 1/P = 1 +E >E$, as expected. 

In the subsequent analysis, we will compare p-variables and e-variables
for other hypotheses 
with those in Theorem \ref{th:1}. 
For a concise presentation, we will always write
\begin{align}\label{eq:baseline}
P_0=(1+(X-\mu)_+^2/\sigma^2)^{-1} \mbox{~~~and~~~} E_0=(X-\mu)_+^2/\sigma^2,
\end{align}
which are the p-variable and e-variable in Theorem \ref{th:1},
and note the connection $P_0=(1+E_0)^{-1}$.

We next consider the hypothesis $H_{\rm S}(\mu,\sigma) $  of symmetric distributions.

\begin{theorem}\label{th:2}
A semi-precise p-variable for $H_{\rm S} (\mu,\sigma)  $ is $P=  \min  \{ (2E_0 )^{-1} ,  P_0    \}, $ 
and a precise   e-variable for $H_{\rm S}(\mu,\sigma) $ is  $E=2E_0$.
Precise p-variables do not exist for $H_{\rm S}(\mu,\sigma)$.
\end{theorem}

From Theorem \ref{th:2}, the e-variable  for $H_S(\mu,\sigma^2)$, which we denote by $E_{\rm S}$ is improved by a factor of two from  $E_0$ for $H(\mu,\sigma^2)$ due to the additional assumption of symmetry. On the other hand, the p-variable in Theorem \ref{th:2}, denoted by $P_{\rm S}$, is improved from $P_0$  by taking an extra minimum with $1/E_{\rm S}$.
In the most relevant case that $P_0\le 1/2$, or equivalently, $E_0\ge 1$, indicating some evidence against the null, 
we have $P_{\rm S}=1/E_{\rm S}$.

Next, we will see that the hypothesis $H_{\rm U}(\mu,\sigma) $  of unimodal distributions
admits the same precise e-variable but a quite improved p-variable, compared to  $P_0$ and $E_0$.
This class includes, for instance, the commonly used gamma, beta, and log-normal distributions. 

\begin{theorem}\label{th:3}
A precise p-variable for $H_{\rm U} (\mu,\sigma)  $ is   $$P=\max\left\{ \frac 4 9P_0, \frac{4P_0 -1}{3}  \right\},$$
and a precise   e-variable for $H_{\rm U}(\mu,\sigma) $ is  $E=  E_0$.
\end{theorem}

We denote the p-variable in Theorem \ref{th:3} by $P_{\rm U}$ and the e-variable by $E_{\rm U}$.
If $P_0$ is smaller than $3/8$, corresponding to $(X-\mu)/\sigma >(5/3)^{1/2}$, then $P_{\rm U}=4P_0/9$; that is, the unimodality assumption reduces the p-variable by a multiplicative factor of $4/9$ compared to $H(\mu,\sigma)$. On the other hand, the e-variable $E_{\rm U}$ does not get improved at all compared to $E_0$.

The proof of Theorem \ref{th:3}, in particular on the factor of $4/9$ for the p-variable, is based on Theorem 1 of \cite{BKV20},  which gives 
 $$
\sup_{Q\in H_{\rm U}(0,1)} T^Q_X(1-\alpha )  =  \max\left\{ \left(\frac{4 - 9\alpha}{9\alpha }\right)^{1/2}  ,   \left(\frac{3-3\alpha}{1+3\alpha}\right)^{1/2}\right\}  \mbox{~~~for $\alpha \in (0,1)$},
 $$
 and applying Lemma \ref{lem:2}  by inverting  of the above curve  as a function of  $\alpha$.

Finally, we consider the hypothesis $H_{\rm US}(\mu,\sigma) $  of unimodal-symmetric distributions. This class includes, for instance, the popular normal, t-, and Laplace distributions.
To construct a semi-precise p-variable for this hypothesis, we will use the following lemma of quantile bounds within $H_{\rm US}(\mu,\sigma)$, which may be of independent interest.
In what follows, $\id$ is the indicator function; that is, $\id_A(x)=1 $ if $x\in A$ and $\id_A(x)=0$ otherwise. 

\begin{lemma}\label{lem:3}
For $\alpha \in (0,1)$, it holds that 
$$
\sup_{Q\in H_{\rm US}(0,1)}T^Q_X(1-\alpha) =    \left(\frac{2}{9\alpha} \right)^{1/2}\id_{(0,1/6]}(\alpha)  +    3^{1/2}{(1-2\alpha)}\id_{(1/6,1/2]} (\alpha).
$$
\end{lemma}
The general formula for $H_{\rm US}(\mu,\sigma)$ 
can be easily obtained  from Lemma \ref{lem:3} via $$\sup_{Q\in H_{\rm US}(\mu,\sigma)}T^Q_X(1-\alpha) =\mu+\sigma \sup_{Q\in H_{\rm US}(0,1)}T^Q_X(1-\alpha) .$$

\begin{theorem}\label{th:4}
A semi-precise p-variable for $H_{\rm US} (\mu,\sigma)  $ is $$P= \frac{2}{9E_0} \id_{[4/3,\infty)}(E_0) +  \frac{3-(3 E_0)^{1/2}}{6}\id_{ (0,4/3)}(E_0) + \id_{\{0\}}(E_0). $$
and a precise   e-variable for $H_{\rm US}(\mu,\sigma) $ is  $E=2 E_0$. 
Precise p-variables do not exist for $H_{\rm US}(\mu,\sigma)$.
\end{theorem} 

The proof of Theorem \ref{th:4}   relies on Lemma \ref{lem:3},  which is a new technical result. The value $2/9$ appeared earlier  in Table 1 of \cite{LSWY18} for $\alpha \le 1/6$,  a  result weaker than Lemma \ref{lem:3}. 

We denote the p-variable obtained from Theorem \ref{th:4} by $P_{\rm US}$ and the e-variable by $E_{\rm US}$. 
 One may check that $P_{\rm US}$ is   smaller than both $P_{\rm U}$ and $P_{\rm S}$ unless $X\le \mu$, in which case they are equal to $1$. 
 For $(X-\mu)/\sigma\ge (5/3)^{1/2}$, or equivalently, $P_0\le 3/8$, we have the following simple relation:
 $$P_{\rm S} = \frac{P_0}{2(1-P_0)},~~ {P_{\rm U}}= \frac{4}{9}P_0,~~\mbox{and}~~P_{\rm US}= \frac{2P_0}{9(1-P_0)},
 $$
implying the order $P_0 >  P_{\rm S} > P_{\rm U} > P_{\rm US}$ unless $P_0=0$. For instance, if  we observe $(X-\mu)/\sigma=3$, then the p-values are
$$
P_0=\frac 1{10}=0.1,~~P_{\rm S} = \frac{1}{18}\approx 0.056,~~ {P_{\rm U}}= \frac{2}{45}\approx 0.044,~~\mbox{and}~~P_{\rm US}= \frac{2}{81}\approx 0.025.
$$
On the other hand,   the corresponding e-values are
$$
E_0=9,~~E_{\rm S}=18,~~E_{\rm U}=9,~~\mbox{and}~~E_{\rm US}=18.
$$
For a comparison, if we are testing the simple parametric hypothesis $\mathrm{N}(0,1)$ against $\mathrm{N}(3,1)$ with one observation $X=3$, then the corresponding Neyman-Pearson p-value is $0.00135$ 
and the corresponding  likelihood ratio e-value  is $90.02$.
This is not surprising as generally p-values and e-values built for composite hypotheses are more conservative than those for simple hypotheses based on the same data.


We summarize our construction formulas for p-variables and e-variables in Table \ref{tab:1} by breaking them down using ranges of $X$.
To obtain the formulas for a general $(\mu,\sigma)$ other than $(0,1)$, it suffices to replace $X$  in  Table \ref{tab:1} by $(X-\mu)/\sigma$.

\begin{table}[t]
\caption{Formulas for p-variables and e-variables}
    \centering
    \renewcommand{\arraystretch}{1.2}
    \tabcolsep=0.5cm
    \resizebox{12cm}{!}{
    \begin{tabular}{c|cl|c}
    \hline
   Hypothesis & \multicolumn{2}{c|}{p-variable} & e-variable \\\hline
  $H(0,1)$  & \multicolumn{2}{c|}{$(1+X_+^2)^{-1}$} & $X_+^2$\\\hline
     \multirow{2}{*}{$H_{\rm S}(0,1)$}   &$    \frac 12   X  ^{-2}  $ & if $X\ge 1$ &    \multirow{2}{*}{ $2X_+^2$}\\ &$   (1+X_+^2)^{-1}     $ & if $X<1$  &  \\ \hline
          \multirow{2}{*}{$H_{\rm U}(0,1)$}   & 
          $   \frac 4 9(1+X ^2)^{-1} $  & if $X\ge (5/3)^{1/2}$ &     \multirow{2}{*}{$X_+^2$}   \\ 
          & $  \frac{4}{3} (1+X_+^2)^{-1} -\frac{1}{3}   $
          & if $X< (5/3)^{1/2}$  &  \\ \hline
         \multirow{3}{*}{  $H_{\rm US}(0,1)$   }  &$\frac{2}{9 } X  ^{-2} $ &if $X\ge (4/3)^{1/2} $ &   \multirow{3}{*}{$2X_+^2$} \\
         &$   \frac 12 - \frac{3^{1/2} }{6}X $ 
         &if $0<X<(4/3)^{1/2} $&     \\
         &$1 $ & if $X\le 0$
         &    \\\hline
    \end{tabular}
    }
    \label{tab:1}
\end{table} 

We conclude the section by making a few technical remarks on the obtained results.

 First, all results holds true if the conditions $\E^Q[X] \le \mu$ and $\var^Q(X) \le \sigma^2$ in each hypothesis is replaced by $\E^Q[X] = \mu$ and $\var^Q(X) = \sigma^2$, respectively. 
Such modifications narrow the hypotheses and 
hence all validity statements hold. 
The precision statements  can be checked with similar arguments to our proofs, and we omit them.
Therefore, knowing $\var^Q(X) = \sigma^2$ on top of $\var^Q(X) \le \sigma^2$, or $\E^Q[X] = \mu$ on top of $\E^Q[X] \le \mu$, does not lead to more powerful one-sided p-variables or e-variables.   

Second, admissibility of the proposed p-variables and e-variables needs future research. 
For e-variables,  admissibility is not difficult to establish, but the picture is different for p-variables.
By Lemma \ref{lem:1},  there always exists a smallest p-variable. 
It remains unclear whether the p-variables we obtained in Theorems \ref{th:1}-\ref{th:4} are the smallest ones for the four hypotheses, respectively.   

Third, for any hypothesis $H$, we can define a  function $g: \alpha \mapsto \sup_{Q\in H} T_X^Q(1-\alpha)$. 
If $g$ is strictly decreasing on $(0,1)$, as in the case of $H(\mu,\sigma)$ and $H_{\rm U}(\mu,\sigma)$, then choosing $f=g^{-1}$ 
yields a precise p-variable $f(X)$. 
For  $H$ being $H_{\rm S}(\mu,\sigma)$ and $H_{\rm US}(\mu,\sigma)$,
   $g$  is flat on $[1/2,1)$,
making   it impossible to find a decreasing $f$ such that  $\inf_{Q\in H} T_{f(X)}^Q(\alpha)=\alpha$ for all $\alpha\in (0,1)$. 
 
\section{Testing the null hypotheses}\label{sec:multi}

\subsection{Constructing e-processes}
\label{sec:41}

We next build tests based on e-values and p-values in Section \ref{sec:single}. 
Section \ref{sec:41} describes the main   methodology based on e-processes for the one-sided testing problem; 
Section \ref{sec:42} describes a few other methods using our results in Section \ref{sec:single}; and 
Section \ref{sec:two-side} discusses the two-sided testing problem on the mean with given variance.

Let
$\mu \in \R$ and $\sigma >0$. We consider the following hypotheses by keeping the same notation as in Section \ref{sec:single}: 
 \begin{align*}
H(\mu,\sigma)&=\left \{Q :\E^Q[X_i|\mathcal F_{i-1}]\le \mu \mbox{~and~}\var^Q(X_i|\mathcal F_{i-1}) \le \sigma^2 \mbox{ for  $i\in [n]$}\right\};
 \\ 
H_{\rm S}(\mu,\sigma)&=\{Q\in H(\mu,\sigma): \mbox{$X_i|\mathcal F_{i-1}$ is symmetrically distributed for  $i\in [n]$}\};
 \\ 
H_{\rm U}(\mu,\sigma)&=\{Q\in H(\mu,\sigma): \mbox{$X_i|\mathcal F_{i-1}$ is unimodally distributed for  $i\in [n]$}\};
 \\
H_{\rm US}(\mu,\sigma)&=H_{\rm U} (\mu,\sigma)\cap H_{\rm S}(\mu,\sigma).
 \end{align*}
 Recall that it is without loss of generality to consider $\mu$ and $\sigma^2$ as constants.
We can also test the hypotheses where some data are symmetric or unimodal and some are not, because we will build e-values from each of them separately.  For simplicity, we only list the above four representative cases. 
Using a similar formulation, the hypothesis in \cite{WR24} is  $$H_{\mathrm{WSR}}(\mu)=\{Q\in H(\mu,1): \mbox{$X_i|\mathcal F_{i-1}$ is supported in $[0,1]$ almost surely for  $i\in [n]$}\}.  $$
In the above formulation, the choice of $\sigma=1 $ is simply to remove the variance constraint; see Remark \ref{rem:R1-WR}.
%
%

There are several simple ways to use results in Section \ref{sec:single} to construct an e-variable or p-variable for the above hypotheses; some of these methods are more useful than the others. 
In general, we can compute an e-variable $E_i$ or p-variable $P_i$ based on $X_i$ for $i\in [n]$  using  Theorems \ref{th:1}-\ref{th:4}, and then combine them. 

Our main proposal is to use e-processes. 
An e-process $M=(M_t)_{t\in [n]}$ can be constructed using 
\begin{align}\label{eq:e-process}
M_t=\prod_{i=1}^t (1-\lambda_i+ \lambda _i E_i),
\end{align}
where $\lambda_i $ is  $\mathcal F_{i-1}$-measurable and takes values in $[0,1)$.
This idea is the main methodology behind game-theoretic statistics; see  \cite{S21}, \cite{SV19}, and \citet[Proposition 3]{WR24}.
 It has been used by \cite{WR24} for testing the mean and \cite{WWZ22} for testing risk measures. 
 To find good choices of $\lambda =(\lambda_i)_{i\in [n]}$ is a non-trivial task.
We propose to specify $\lambda$ in two different ways. 
\begin{enumerate}[(a)]
\item \textbf{E-mixture method}: 
We first take several $\lambda_i=\lambda \in [0,1)$, which is a constant for each $i\in [n]$, and then average the resulting e-processes from \eqref{eq:e-process} over these choices to get an e-process.
An uninformative  choice of the values of $\lambda$ 
may be some points in $[0,0.2]$.
We avoid choosing $\lambda$ close to $1$ because our e-value may take the value $0$ with substantial probability, leading a small value of $\E^Q[\log (1-\lambda + \lambda E)]$. This quantity measures the growth rate of an e-process; see \cite{GDK20} and \cite{WR24}.
In our simulation and empirical studies, we average over $\lambda =0.01\times \{1,\dots,20\}$.

\item \textbf{E-GREE method}: In the GREE (growth-rate for empirical e-statistics) method of \cite{WWZ22} for $\lambda_i$, $i\in [n]$ in \eqref{eq:e-process}, $\lambda_i$ is determined by solving the following optimization problem:
\begin{equation}
\label{eq:gree}
\lambda_i = \left(\underset{{\lambda \in [0, 1)}}{\arg\max} \frac{1}{i -1} \sum_{j= 1}^{i - 1}\log(1 - \lambda + \lambda E_j) \right)\wedge \frac 12 .
\end{equation}
To simplify the maximization in \eqref{eq:gree}, a fast and approximate solution can be obtained using Taylor expansion as in \cite{WR24}. This leads to the following simple formula
\begin{equation}
\label{eq:agree}
\lambda_i =  \left( \frac{\sum^{i-1}_{j = 1} (E_j - 1)}{\sum^{i-1}_{j = 1} (E_j - 1)^2} \right)_+\wedge \frac 12.
\end{equation}
We will use \eqref{eq:agree} for all e-GREE related calculations for the following results. Our unreported simulation   suggests that using \eqref{eq:gree} and using \eqref{eq:agree} yield very similar results.
\end{enumerate}
When the hypothesis to test is $H_{\rm WRS}(\mu)$, the e-GREE method reduces to the method of \cite{WR24}; see Section \ref{sec:52}.
An optimization procedure related to \eqref{eq:gree} is studied by \cite{KTT11}.

For either the e-GREE or the e-mixture method, we fix $\alpha\in (0,1)$ and reject the null hypothesis if the e-process $M$ goes beyond $1/\alpha$, that is, when $M_t \ge 1/\alpha$ for the first time. 
The Type-I error control is guaranteed by Ville's inequality (\cite{V39}) as 
$\p(\sup_{t\in [n]}M_t\ge 1/\alpha)\le \alpha$, because any e-process is almost surely upper bounded by nonnegative supermartingales with initial value one; see \cite{RRLK20}.

The result below clarifies   consistency of the e-GREE method in the most idealistic setting. 

\begin{proposition}\label{prop:consistency}
Suppose that data are  iid  and   generated from an alternative probability $Q$. The e-GREE method has asymptotic power approaching $1$ as $n\to \infty$, that is, $Q(\sup_{t\in [n]}M_t \ge 1/\alpha) \to 1$ for any $\alpha\in (0,1)$ if and only if $\E^Q[E_1]>1$.   
\end{proposition}

Although Proposition \ref{prop:consistency} requires an iid assumption, this assumption is not needed for consistency in practical situations; a simulation example is in Section \ref{sec:51}.

\subsection{Some other methods}
\label{sec:42}

Below we list some other methods, where we assume that $n$ is finite. 
They generally do not work well as shown by the simulation studies, but nevertheless we list them as they follow from our results in Section \ref{sec:single}, and they are presented only for a comparison.

\begin{enumerate}[(a)]
\item[(c)] \textbf{P-Fisher method}: Construct a p-variable $P$ using the Fisher combination
$$
P =1- \chi_{2n}(-2(\log   P_1 +\dots +\log P_n)),
$$ 
where $\chi_{2n}$ is the cdf of a chi-square distribution with $2n$ degrees of freedom.

\item[(d)] \textbf{P-Simes method}: Construct a p-variable $P$ using the Simes combination; see \cite{S86}, 
$$
P = \min_{i\in [n]} \frac{n }{i } P_{(i)},
$$ 
where $P_{(i)}$ is $i$-th order statistic of $P_1,\dots,P_n$ from the smallest to the largest. 
\end{enumerate} 
Although in general p-Fisher and p-Simes require independence among p-variables,  they are valid in our setting since our p-variables are conditionally valid, and they can be combined as if they are iid; a proof of this is presented in the Supplementary Material.

Then next two methods use all data directly, and requires independence among $X_1,\dots,X_n$. 
A most natural statistic is the sample mean $T=\sum_{i=1}^n X_i/n$. 
Under $H(\mu,\sigma)$,   $T$ has at most mean $\mu$ and variance at most $\sigma^2/n$.
Moreover, symmetry of $T$ follows from symmetry of $X_1,\dots,X_n$.
Nevertheless, $T$ is not necessarily unimodal even if $X_1,\dots,X_n$ are unimodal, and hence unimodality of $T$ cannot be used. 
The following e-variables and p-variables are constructed by directly applying  Theorems \ref{th:1}-\ref{th:4}. 
\begin{enumerate}[(a)]
\item[(e)] \textbf{E-batch method}: An e-variable for   $H(\mu,\sigma)$ or  $H_{\rm U}(\mu,\sigma)$ is 
$$E_0 = n (T-\mu)_+^2/\sigma^2, $$ an e-variable for   $H_{\rm S}(\mu,\sigma)$ or $H_{\rm US}(\mu,\sigma)$ is 
$$E_{\rm S} = 2n (T-\mu)_+^2/\sigma^2.$$

\item[(f)] \textbf{P-batch method}: A p-variable for   $H(\mu,\sigma)$  or $H_{\rm U}(\mu,\sigma)$ is 
$$P_0= (1+E_0)^{-1},$$
a p-variable for  $H_{\rm S}(\mu,\sigma)$ or $H_{\rm US}(\mu,\sigma)$ is 
 $$P_{\rm S} = \min  \{ (2E_0 )^{-1} ,  P_0    \}.$$
 \end{enumerate}
 All methods described in this section have Type-I error control under the null hypothesis and with finite sample (with methods (e) and (f) additionally requiring independence)
 without requiring  that the data are identically distributed.

\subsection{Two-sided e-values testing the mean given variance}
\label{sec:two-side}

We briefly discuss the two-sided mean testing problem, where the main hypothesis $H(\mu^L,\mu^U, \sigma)$ to test is 
$$
\left\{Q :\E^Q[X_i|\mathcal F_{i-1}]\in [\mu^L,\mu^U] \mbox{~and~}\var^Q(X_i|\mathcal F_{i-1}) \le \sigma^2 \mbox{ for  $i\in [n]$}\right\},
$$
where $\mu^L\le \mu^U$ are constants. The case $\mu^L=\mu^U$ corresponds to testing whether the mean is equal to a precise value.

Our methodology can be easily adapted to test this hypothesis. 
First, we note that the e-variable $E$ given by 
\begin{equation}
    \label{eq:e-twoside}E = \frac{(X-\mu^U)_+^2 +(X-\mu^L)_-^2 }{\sigma^2} ,
\end{equation} 
is a precise e-variable for $H(\mu^L,\mu^U, \sigma)$ formulated on a single observation $X$. 
To see this, it suffices to note that for $Q\in  H(\mu^L,\mu^U, \sigma)$, 
\begin{align*}
\E^Q[E] &= \E^Q\left[ \frac{(X-\mu^U)_+^2 +(X-\mu^L)_-^2 }{\sigma^2} 
\right]\\
&\le 
\E^Q\left[ \frac{(X-\E^Q[X])_+^2 +(X-\E^Q[X])_-^2 }{\sigma^2} 
\right]
=\frac{\var^Q(X)}{\sigma^2}\le 1.
\end{align*}
The statement on its precision can be verified similarly to Theorem \ref{th:1}.

If $\mu^L=\mu^U=\mu$, then the e-variable in \eqref{eq:e-twoside} is 
$$
E={(X-\mu)^2}/{\sigma^2}.
$$
This e-variable satisfies the property that $\E^Q[E]>1$ if $\E^Q[X]\ne \mu$
and $\var^Q(X)=\sigma^2$; this condition is useful to establish consistency in Proposition \ref{prop:consistency}.

Following the same procedure in Section \ref{sec:41} using \eqref{eq:e-twoside}, we obtain e-processes for the two-sided problem $H(\mu^L,\mu^U, \sigma)$. 
Due to a smaller null hypothesis, 
this e-process is generally more powerful than the one in Section \ref{sec:41} testing the one-sided mean. 

There are special, adversarial scenarios where such two-sided tests may not be powerful. 
For instance, if  data are independent   with $\E[X_i]<\mu$ and $\E[X_j]>\mu$ appearing in an alternating sequence; this forms a dataset that looks like iid data with mean $\mu$, thus very difficult to detect. The same challenge   exists for other methods based on e-processes, such as that of \cite{WR24}.
 \begin{remark}
Under the additional information of symmetry, the e-variable in \eqref{eq:e-twoside} 
 can be used, but  it  cannot be multiplied by two as in Theorem \ref{th:2}.
 In this case, an alternative way to take advantage of symmetry is to   build two e-processes in Section \ref{sec:41}:  one to test $\E[X_i|\mathcal F_{i-1}] \le \mu^U$
and another one to test $\E[-X_i|\mathcal F_{i-1}] \le - \mu^L$. 
Taking the average of these two e-processes yields a valid e-process for the null hypothesis. As long as one of the two e-processes has good power for the true data generating procedure,  the average e-process has good power. 
\end{remark}

\subsection{Power of the e-values with fixed mean and growing variance}
\label{sec:r1-add-1}

In this section, we analyze   the power of the e-variables.
For a given e-variable $E$, its e-power, using the terminology of \cite{VW24}, for an alternative probability $Q$ is defined as $\E^Q[\log E]$; see \cite{S21} and \cite{GDK20} for using this quantity as a notion of power. 
Certainly, the power depends on the specific alternative $Q$. 
We are particularly interested in how the e-power changes as the variance  in the alternative hypothesis grows. 

For this purpose, we consider a simplistic, yet representative setting, where a class of simple alternatives $(Q_\sigma)_{\sigma >  1}$ is  indexed by $\sigma > 1$, 
such that our data point $X$ under $Q_\sigma$ is distributed as $\sigma Z$, where $Z$ has a fixed distribution with mean $0$ and variance $1$ satisfying the null hypothesis, which can be  one of $H(0,1)$, $H_{\rm S}(0,1)$, $H_{\rm U}(0,1)$ and $H_{\rm US}(0,1)$.
Note that in this setting, the mean of the data is always $0$, and only its variance grows  under the alternative.
We denote by $Q_0$ a null probability.
Below, we will show that  the e-power of each  e-variable grows at a rate of $\log \sigma$ as the alternative variance $\sigma^2$ grows, 
regardless of the distribution of $Z$.

Let $E$ be the e-variable computed based on $X$ as in  Section \ref{sec:single}.  
Due to the construction of the e-process $M$ in \eqref{eq:e-process},
the e-power of relevance is defined as
$$
\Pi^{Q_\sigma}= \sup_{\lambda \in [0,1]} \E^{Q_\sigma}[\log ( 1-\lambda + \lambda E)] = \sup_{\lambda \in [0,1]} \E^{Q_0} [\log ( 1-\lambda + \lambda \sigma^2 E)],
$$
that is, the best-achievable e-power of  in each multiplicative term in the e-process $M$. 
\begin{proposition}\label{prop:power}
Suppose $p:=Q_0(E\ge 1)>0$. For $\sigma>1$,
\begin{align}\label{eq:rw-r1-1}
( 2 p   \log \sigma -\log 2 )_+  \le \Pi^{Q_\sigma} \le  2 \log \sigma.
\end{align}
Moreover, $ 0\le \Pi^{Q_\sigma} -\Pi^{Q_\delta} \le2( \log \sigma -\log \delta) $ for $\sigma>\delta>1$.
\end{proposition}

Proposition \ref{prop:power} suggests that the growth rate of the e-process $M$ 
is roughly a constant times $\log \sigma$ when the alternative  variance $\sigma^2$ is larger than $1$.
An additional negative term $-\log2$  in \eqref{eq:rw-r1-1} is not surprising, because our conditions 
do not guarantee $\Pi^{Q_\sigma}>0$ for $\sigma$ very close to $1$.
Below, we give an example to illustrate the sharpness of bounds in \eqref{eq:rw-r1-1}.

\begin{example}
Suppose that $Q_0(E=0)=Q_0(E=2)=1/2$. We can compute 
\begin{align*}
 \Pi^{Q_\sigma} 
 =  \sup_{\lambda \in [0,1]} \frac 12 \left( \log ( 1-\lambda )+  \log ( 1 +  \lambda (2\sigma^2-1))  \right)
 = \frac 12 \log \frac{\sigma^4}{2\sigma^2-1}.
 \end{align*}
 It is clear that $ \Pi^{Q_\sigma} $ is approximately equivalent to $ \log \sigma$ for  large $\sigma$, corresponding to the left side of \eqref{eq:rw-r1-1} with $p=1/2$.
\end{example}

 

\section{Simulation studies}\label{sec:simulation}

\subsection{A comparison of different e-combining methods}
\label{sec:51}
In this section, we conduct simulation studies for the non-parametric hypotheses in Section \ref{sec:multi}. We  set  $\mu=0$ and $\sigma = 1$ without loss of generality. 
  
We first concentrate on the null hypothesis $H(0,1)$, as the other four cases are similar. For all the methods stated in Section \ref{sec:multi}, we do not make the assumption that the data are identically distributed. Thus, we generate a sample of $n$ independent data points, although independence is not needed for methods (a)-(d), 
alternating from two different distributions:  $X_1,X_3,\dots,$ follow a normal distribution,
and $X_2,X_4,\dots,$ follow a Laplace distribution, with the same mean $\nu$ and the same variance $\eta^2$.\footnote{The assumption that the two distributions have the same mean and variance is not necessary when evaluating the power of the methods. We assume this only for simplicity.}  
We denote this data generating process as $\mathrm{NL}(\nu,\eta^2)$ with the null parameters being $(\nu,\eta^2)=(0,1)$.  
We consider two alternatives: (1) Data generated from $\mathrm{NL}(0,\eta^2)$ where $\eta>1$;  (2) Data generated from $\mathrm{NL}(\nu,1)$ where $\nu>0$.  
In our setting, the tester does not know  the alternating data generating mechanism. For each alternative model, we compute the rejection rate over 1000 runs using the thresholds of $E\ge 1/\alpha$ and $P\le \alpha$, with $\alpha=0.05$, for e-values and p-values, respectively.

 For the e-mixture method, we experiment by averaging $\lambda$ in the interval $[0.01, 0.20]$ with step size $0.01$. The e-GREE method is similar to the e-mixture method, except that $\lambda_i$ is dynamically updated with different $i\in [n]$ using the formula \eqref{eq:gree}. 
 
\begin{figure}[!t]
\centering
\includegraphics[width=\textwidth]{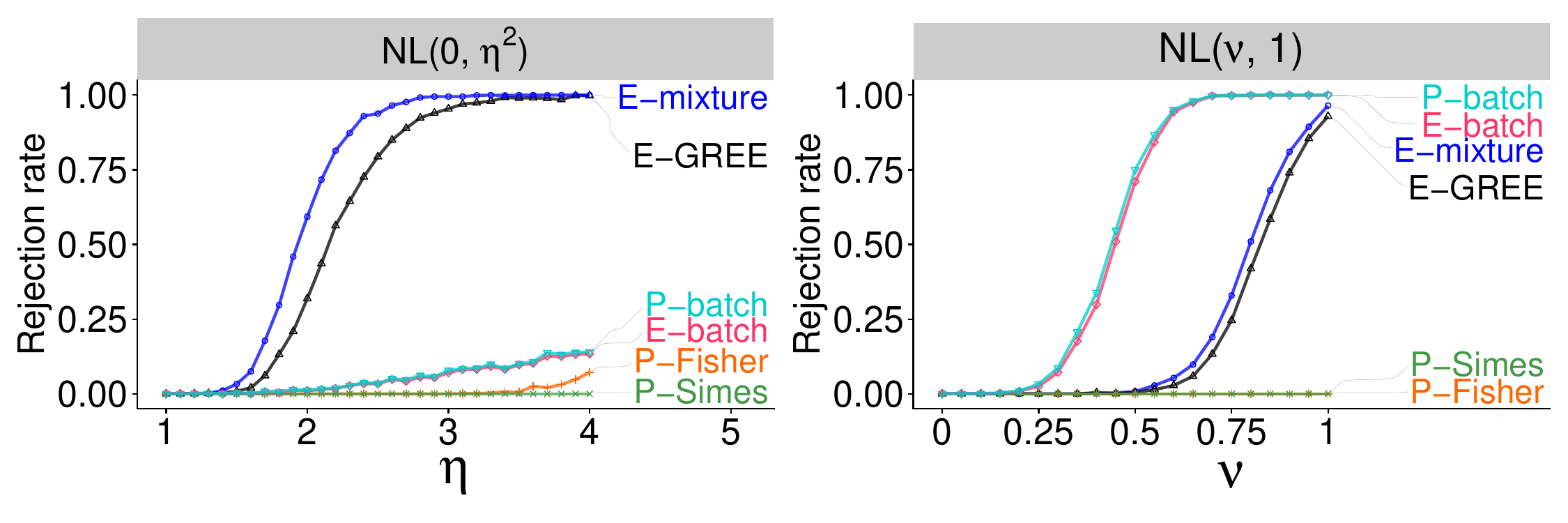}
\caption{Rejection rates for all methods for testing $H(0,1)$ with sample size $n = 100$ over 1000 runs using the threshold $ 20$.} 
\label{fig:rej}
\end{figure}

Figure \ref{fig:rej} shows the rejection rates for all methods with data generated from $\mathrm{NL}(0,\eta^2)$ for $\eta\in [1,4]$, and from $\mathrm{NL}(\nu, 1)$ for $\nu\in [0,1]$. For alternative model $\mathrm{NL}(0,\eta^2)$, we see that the e-mixture   and the e-GREE methods outperform the other methods, with the e-mixture method being the most powerful. For $\eta<1.5$, the rejection rates of all methods are very low, making it challenging to distinguish their efficiency. As $\eta > 1.5$, both the e-mixture method and the e-GREE method exhibit significantly higher rejection rates compared to other methods, demonstrating their effectiveness in testing $H(0,1)$. The other four methods have almost no power.  For alternative model $\mathrm{NL}(\nu, 1)$, we observe that e-batch method and the p-batch method show significant high rejection rates, since they are quite sensitive to  the sample mean. Recall that these methods rely on independence, so the central limit theorem kicks in.

Among all methods, only the e-process based methods satisfy anytime-validity, that is, decision can be made at any stopping time when data arrive sequentially.  This situation is common in financial applications, where realized losses accumulate over time; see the empirical study in Section \ref{sec:emp}.

The testing procedures for $H_{\rm S}$, $H_{\rm U}$ and $H_{\rm US}$ are the same as for testing $H$. We generate 100 data points from $\mathrm{NL}(0.5,2)$ and calculated the rejection rates for testing $H_{\rm S}$, $H_{\rm U}$ and $H_{\rm US}$ with null hypotheses $\mu = 0$ and $\sigma = 1$.  Table \ref{tab:2} displays the rejection rates for all hypotheses. It is clear that the extra information of symmetry improves the power.

\begin{table}[htbp]
\caption{Rejection rates of testing $H$, $H_{\rm S}$, $H_{\rm U}$ and $H_{\rm US}$ with $n=100$ data generated from the model $\mathrm{NL}(0.5,2)$.}
\centering
\renewcommand{\arraystretch}{0.8}
\tabcolsep=0.4cm
\resizebox{12cm}{!}{
\begin{tabular}{ccccccc}
\toprule
&E-mixture & E-GREE & P-Fisher &P-Simes &E-batch & P-batch\\
\midrule
$H$ & 0.419 & 0.315 &  0.000  &   0&  0.639 & 0.664\\
$H_{\rm S}$  & 0.998 & 0.882  &0.000    & 0  &0.900  &0.900\\
$H_{\rm U}$ &0.419 & 0.315 & 0.006   &  0  &0.639  &0.664\\
$H_{\rm US}$ &0.998  &0.882 & 0.763   &  0 & 0.900  &0.900\\
\bottomrule
\end{tabular}
}
\label{tab:2}
\end{table}

\subsection{A comparison with the GRAPA method}
\label{sec:52}
Recall that our model can also be interpreted as testing the mean under the knowledge of an upper bound on the variance. This allows us to compare our testing approach with the GRAPA (Growth Rate Adaptive to the Particular Alternative) method proposed by \cite{WR24}. GRAPA is similar to the e-GREE method discussed in Section \ref{sec:multi}, but it requires the random variable to be bounded. The e-process $(M_t)_{t \in [n]}$ for the GRAPA method is constructed as follows:
\begin{align}\label{eq:grapa}
M_t=\prod_{i=1}^t (1 + \lambda_i( X_i - \mu)),
\end{align}
where $\mu$ is the conditional mean being tested and $\lambda_i$ is $\mathcal F_{i-1}$-measurable and takes value in $(-1/(1-\mu),1/\mu)$. 
It is clear that $1 + \lambda_i( X_i - \mu)$ is an e-variable for each $i\in [n]$.  Thus, maximizing the growth of \eqref{eq:grapa} is similar to \eqref{eq:gree}, where $\lambda_i$ is determined by solving the following optimization problem:
\begin{equation}
\label{eq:egrapa}
\lambda_i = \underset{{\lambda \in [-c/(1-\mu), c/\mu]}}{\arg\max} \frac{1}{i -1} \sum_{j= 1}^{i - 1}\log(1 + \lambda( X_i-\mu)),
\end{equation}
where $c\in (0,1]$ is fixed. 
For faster computation in the context of confidence sequences, 
\cite{WR24} also offered an alternative way to obtain $\lambda_i$, which they called approximate GRAPA method, and $\lambda_i$ is determined by 
\begin{equation}
\label{eq:agrapa}
\lambda_i = - \frac{c}{1-\mu} \vee \frac{\widehat{\mu}_{i-1} - \mu}{\widehat{\sigma}^2_{i-1} + (\widehat{\mu}_{i-1}-\mu)^2} \wedge \frac{c}{\mu},
\end{equation}
where $\widehat{\mu}_{i}$ and $\widehat{\sigma}^2_{i}$ are empirical mean and variance of the observations $X_1 \cdots, X_i$.
From \eqref{eq:agrapa}, it is clear that the GRAPA method is able to use the sample variance information adaptively.
In particular, our e-GREE method in  \eqref{eq:agree} 
is  adaptive to the empirical variance of  the e-values. In the simulation results, we use \eqref{eq:egrapa} and choose $c = 1/2$.


We compare  five  methods for testing the mean under various conditions:
\begin{itemize}
    \item[(a)] GRAPA: The GRAPA method with a bounded support $[0,1]$.
    \item[(b)] E-GREE: The e-GREE method with the variance upper bound $\sigma^2$.
    \item[(c)] E-mixture: The e-mixture method  with the variance upper bound $\sigma^2$. 
  \item[(d)] E-GREE-2s: The two-sided e-GREE method with the variance upper bound $\sigma^2$.
    \item[(e)] E-mixture-2s: The two-sided e-mixture method with the variance upper bound $\sigma^2$. 
\end{itemize} 
 We note that GRAPA is designed as a two-sided test, although it can easily be adjusted by restricting $\lambda_i$ in \eqref{eq:grapa} to be non-negative.
\begin{remark} We could also implement the e-GREE and e-mixture methods without an upper bounded variance but using the bounded support, as described in Remark \ref{rem:variance}. Although these methods are valid, they have poor power in our setting, because their assumption is strictly weaker than both bounded variance and bounded support. We omit their results. 
\end{remark}

 We set $\mu = 0.35$ and apply both one-sided and two-sided tests on the same dataset. We generate a sample consisting of $n$ independent data points from a beta distribution, denoted by $\mathrm{Beta}(\nu,\sigma^2)$, where $\nu$ and $\sigma^2$ represent the mean and variance of the beta distribution.\footnote{None of the methods requires that the data follow identical distributions; we use a single distribution just for simplicity.} Here, we use $\nu$ and $\sigma^2$ instead of the standard beta parameters $\alpha$ and $\beta$ for the sake of convenience.  Note that the parameters $\alpha$ and $\beta$ can be easily recovered based on given mean $\nu$ and variance $\sigma^2$: $\alpha = \nu(\nu - \nu^2 - \sigma^2)/ \sigma^2$ and $\beta = (\nu^2 + \sigma^2 - \nu)(\nu-1)/ \sigma^2$. Since the beta distribution has a  bounded support $[0,1]$,  we can make meaningful comparisons between the GRAPA method and the e-GREE and e-mixture methods.



We first compare the rejection rates, using a threshold of $ 20$ over 1000 runs, for all methods mentioned above under different $\nu$ with fixed $\sigma^2$. We consider  $\nu\ge 0.35$ and $\sigma=0.05$, $\sigma=0.1$ and $\sigma=0.3$.   We use 20 data points for each run.

\begin{figure}[!t]
\centering
\includegraphics[width=\textwidth]{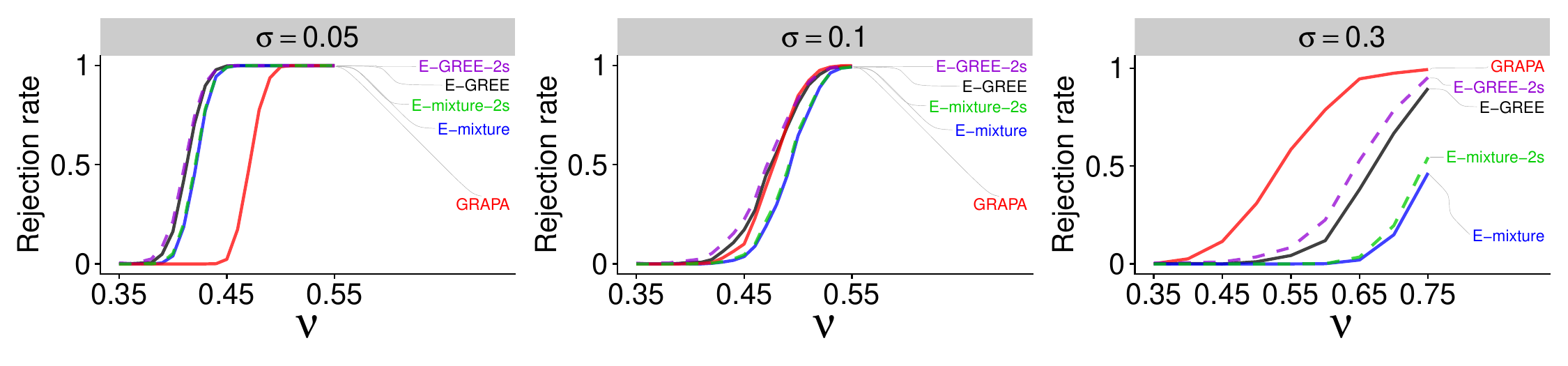}
\caption{Rejection rates for the GRAPA, the e-GREE, the e-mixture and the two-sided e-GREE-2s and the e-mixture-2s methods over 1000 runs using the threshold  $ 20$ and $\mu = 0.35$. Data are generated from $\mathrm{Beta}(\nu,\sigma^2)$ with sample size $n = 20$, where $\nu \ge 0.35$ and $\sigma\in\{0.05,0.1,0.3\}$.}
\label{fig:2full}
\end{figure}

Figure \ref{fig:2full} shows the performance of the three methods. First, the e-GREE method is always better than the e-mixture method.  Second, the two-sided versions of both the e-GREE and e-mixture methods show a slight improvement over their respective one-sided methods, as expected. Third, in case $\sigma=0.05$ and $\sigma=0.1$, the e-GREE method outperforms the GRAPA method; in case $\sigma=0.3$, the GRAPA method demonstrates superior performance compared to the other methods. This is intuitive, because the variance information is less useful for larger $\sigma$; recall that for any distribution supported in $[0,1]$ with mean $\mu\le 0.35$,  the maximum possible variance is $0.2275$, and $\sigma \approx 0.477$.

\begin{figure}[!t]
\centering
\includegraphics[width=\textwidth]{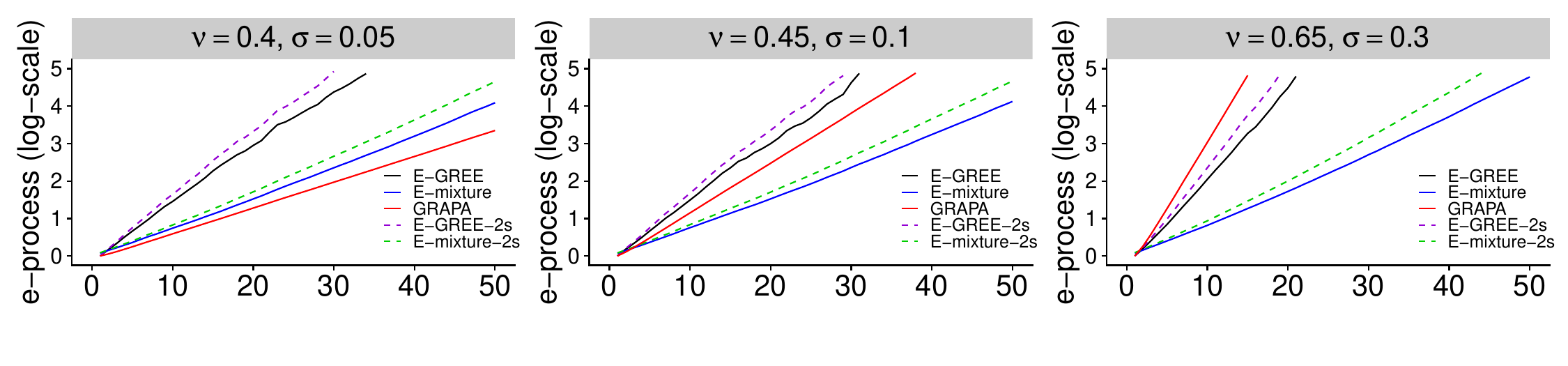}
\caption{Average logarithmic e-processes for the GRAPA, the e-GREE, the e-mixture and the two-sided e-GREE-2s and the e-mixture-2s methods  with varying sample size and $\mu=0.35$. Data are generated from $\mathrm{Beta}(\nu, \sigma^2)$ where $\sigma \in \{0.05,0.1,0.3\}$  and $\nu =\mu + \sigma$.}   
\label{fig:3full}
\end{figure}

Figure \ref{fig:3full} shows the average logarithmic e-processes for $n$ up to $50$ by using $\nu = \mu + \sigma$ for each alternative model. The relative rankings of these methods are consistent with their rejection rates, with e-GREE performing the best when $\sigma$ is relatively small. 
 
From the simulation results, our general recommendation is to use e-GREE to construct the e-process when the variance to be tested is relatively small, and to use GRAPA when the variance to be tested is relatively large compared to the bounded support. 
%

\subsection{A comparison with exponential test supermartingale}

Next, we compare our methods with the exponential test supermartingale methods that directly construct e-processes, rather than using a betting strategy to combine sequential e-variables. 

\cite{WR23} extends the idea from \cite{C12} to construct a non-negative test supermartingale called the Catoni supermartingale to test mean and variance in sequential settings. The test supermartingale is constructed as follows:
\begin{equation}
\label{eq:catoni}
    M_t^{\rm C} = \prod^{t}_{i = 1}\exp\left( \phi(\lambda_i(X_i - \mu)) - \frac{\lambda_i^2 \sigma^2}{2} \right),
    \end{equation} 
where $\phi$ is the influence function and $(\lambda_i)_{i \in [n]}$ is any predictable process. Following the  recommendation of \cite{WR23}, we choose the influence function
$$\phi(x) = 
\begin{cases} 
\log(1 + x + x^2/2), & \text{if } x \geq 0; \\
-\log(1 - x + x^2/2), & \text{if } x < 0.
\end{cases} $$
and $(\lambda_i)_{i \in [n]}$ as
\begin{equation}
\label{eq:lambda}
\lambda_i = \left(\frac{2 \log(1/\alpha)}{i(\sigma^2 + \eta_i^2)}\right)^{1/2} \quad \text{where} \quad \eta_i = \left(\frac{2 \sigma^2 \log(1/\alpha)}{i- 2 \log(1/\alpha)}\right)^{1/2}.
\end{equation}

A different approach by  
\cite{HRMS21} is to use a framework for non-parametric confidence sequences based on the concept of exponential supermartingales. They introduce the concept of a ``sub-$\psi$ process'' in  \citet[Definition 1]{HRMS21}.
Informally, a sub-$\psi$ process is
a pair of $\mathcal{F}_t$-adapted processes $(S_t, V_t)$ such that $S_t$ is the zero-mean deviation of the sample sum from its estimand at time $t$ and  $V_t$ and  $\psi$ make the following process 
\begin{equation}
\label{eq:superm}
  M^{\psi}_t = \exp\{ \lambda S_t - \psi(\lambda) V_t \} 
\end{equation}
dominated by a supermartingale for each $\lambda$ in an interval $[0,\lambda_{\max})$. 
This framework allows for testing mean and variance under a wide variety of assumptions, including bounded supports, self-normalized bounds, and symmetric conditions. We refer to \cite[Appendix J, Table 3]{HRMS21} for a collection of commonly used $\psi$ functions and variance processes for $S_t = \sum_{i = 1}^t (X_i - \mu)$ under various assumptions. We choose two special cases for comparison with our methods: the self-normalized bounds test supermartingale, denoted by $M^{\psi, \mathrm{SN}}_t$, and the symmetric condition test supermartingale, denoted by $M^{\psi, \mathrm{sym}}_t$. For $\lambda \in [0, \infty)$, these test supermartingales are constructed as follows:
\begin{equation}
\label{eq:SN}
  M^{\psi, \mathrm{SN}}_t = \prod^{t}_{i = 1}\exp\left( \lambda (X_i - \mu) - \frac{\lambda^2(X_i - \mu)^2 + 2\sigma^2}{6} \right),
\end{equation}
which also appears in \citet[Section 5]{WR23}, and 
\begin{equation}
\label{eq:sym}
  M^{\psi, \mathrm{sym}}_t = \prod^{t}_{i = 1}\exp\left( \lambda (X_i - \mu) - \frac{\lambda^2(X_i - \mu)^2}{2} \right).
\end{equation}
We follow a simple method of choosing $\lambda$ suggested by \citet[Section 3.2]{HRMS21}, that is,  to use the mixture supermartinagle $\int \exp(\lambda S_t - \psi(\lambda) V_t) \, \mathrm{d}\Phi(\lambda)$ by assuming $\lambda \sim \Phi=\mathrm{N}(0,1)$. Now, we further compare the following methods:
\begin{itemize}
    \item[(f)] WR23-Catoni: The Catoni method with the variance upper bound $\sigma^2$.
    \item[(g)] HRMS21-SN: The self-normalized method with the variance upper bound $\sigma^2$.
    \item[(h)] HRMS21-sym: The sub-$\psi$ method with symmetry, but without variance information.
    \item[(i)] E-GREE-sym: The e-GREE method with the variance upper bound $\sigma^2$ and symmetry.
    \item[(j)] E-mixture-sym: The e-mixture method with the variance upper bound $\sigma^2$ and symmetry.
\end{itemize} 

\begin{figure}[!t]
\centering
\includegraphics[width=\textwidth]{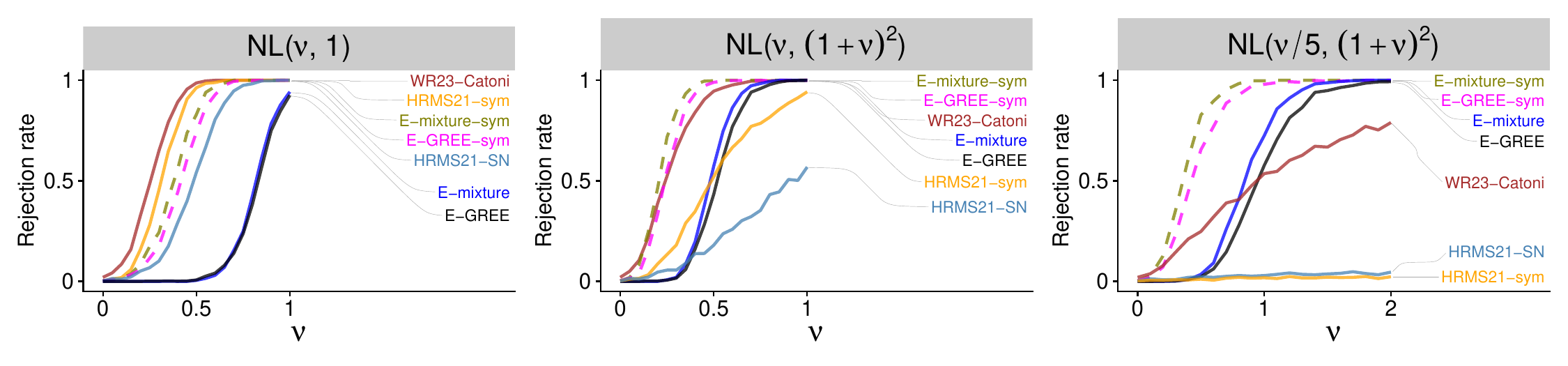}
\caption{Rejection rates for methods (a), (b) and (f)-(j) for testing $H (0, 1)$ with sample size $n = 100$ over 1000 runs using the threshold 20.}
\label{fig:superm}
\end{figure}

We compare above five methods, along with the e-GREE and e-mixture methods that do not utilize symmetric information (methods (a) and (b) described in the previous section), in testing $H(0,1)$. Following the same data generating process as described in Section \ref{sec:52}, we generate $n$ independent data points alternating between the normal and Laplace distributions, denoted by $\mathrm{NL}(\nu, \eta^2)$. Figure \ref{fig:superm} shows rejection rates for above methods with data generated from three cases: $\mathrm{NL}(\nu, 1^2)$ for $\nu \in [0,1]$,     $\mathrm{NL}(\nu, (1+\nu)^2)$ for $\nu \in [0,1]$, and   $\mathrm{NL}(\nu/5, (1+\nu)^2)$ for $\nu \in [0,2]$. 

For $\mathrm{NL}(\nu, 1^2)$, the Catoni method outperforms other methods, while methods utilizing symmetric information generally perform well. For $\mathrm{NL}(\nu, (1+\nu))^2)$, where both the mean and variance of the data generating process change, the power of methods from \cite{HRMS21} reduces. In contrast, the power of our e-value based methods increases, as our construction of e-values is sensitive to the changes to variance. In the last case, $\mathrm{NL}(\nu/5, (1+\nu)^2)$, the impact of changes in mean is small and the variance effect is large, e-value based methods generally outperform others. Although method (h) benefits from not requiring information about variance or even the existence of variance, it demonstrates minimal power when testing mean with varying variance, due to its penalization term $-(X_i-\mu)^2$  in the exponential form of \eqref{eq:SN} and \eqref{eq:sym}.
In summary, our methods are comparatively more powerful when the alternative variance defers from the null. 

\section{Empirical study with financial data}
\label{sec:emp}

In this section, we conduct an empirical study to test the hypothesis $H(\mu,\sigma)$ on the daily losses of financial assets. We aim to calculate the number of trading days required to detect evidence for rejecting the null hypothesis $H(\hat{\mu},\hat{\sigma})$ during the 2007--2008 financial crisis period. Here, $\hat{\mu}$ and $\hat{\sigma}$ represent the sample mean and sample variance estimated from historical data prior to the testing period. 
That is, we are testing whether the historical estimations before the testing period are still valid. 
If the null hypothesis can be rejected at a reasonable thresholds level rather swiftly, this will serve as evidence of the effectiveness of e-process methods and could help investors switch strategies in a timely manner.

We choose 20 stocks from 10 different sectors of the S\&P 500 list with the large market capitalization in each sector.
Moreover, we include two companies with the largest market capitalization from the to-be Real Estate sector.\footnote{Real Estate becomes the 11th sector of S\&P500 in 2016.} 
We first calculate the daily losses for each of the selected stocks from January 1, 2001 to December 31, 2010. The daily losses are expressed by percentage and calculate by $L_t = -(S_{t+1} - S_t)/ S_t$, where $S_t$ is the close price at day $t$.
Note that the positive value represents a loss and negative value represents a gain. 
We could also use the log-loss data instead of the linear loss data, but the difference between the two is minor.   
We use the loss data from  January 1, 2001 to December 31, 2006 to estimate the mean and variance for the null hypothesis. We compute the e-values using both the e-mixture method and the e-GREE method based on the construction of \eqref{eq:e-process} as the daily loss from January 1, 2007 fed into the e-process.  

Following a methodology similar to the simulation study in Section \ref{sec:simulation}, we report the evidence against the null hypothesis when the e-process exceeds thresholds of 2, 5, 10, and 20.\footnote{In accordance with Jeffrey's rule of thumb about e-values (see \cite{J61} and \cite{VW21}), if the e-value falls within the interval of $(10^{1/2}, 10)$, the evidence against the null hypothesis is considered substantial; If the e-value falls within the interval of $(10, 10^{3/2})$, the evidence against the null hypothesis is regarded as strong.} E-values exceeding 5 or 10 provide substantial evidence to reject the null hypothesis, while a threshold of 20 offers strong evidence against the null hypothesis. It is important to note that, although a threshold of 2 may not be substantial enough to reject the null hypothesis, it can still serve as an early warning that the stock's performance may be different from its historical path.

\begin{figure}[t]
\centering
\includegraphics[width=\textwidth]{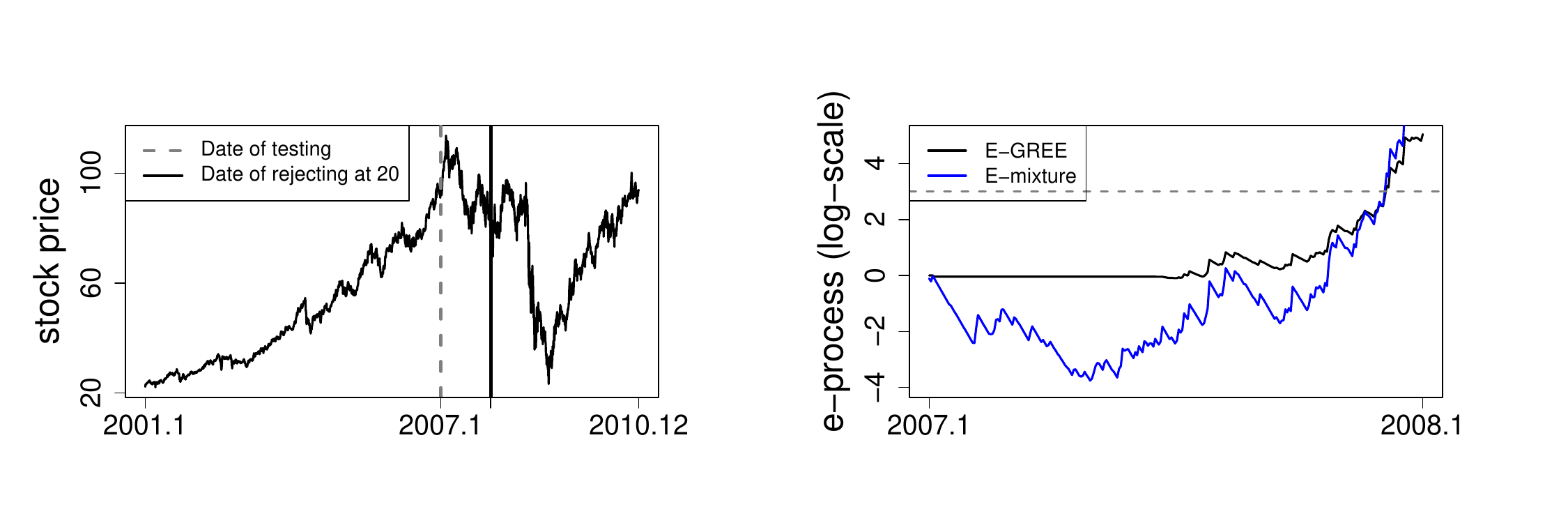}
\caption{Sample path and logarithmic e-process using the e-GREE and the e-mixture methods testing of $H(\hat{\mu}, \hat{\sigma})$ for Simon Property (SPG) stock from January 2007 to January 2008, where $\hat{\mu} = -0.001028$ and $\hat{\sigma} = 0.012123$ are the sample mean and variance estimated from historical data for stock SPG from January 1, 2001 to December 31, 2006.} 
\label{fig:spg}
\end{figure}

To illustrate the e-process detection procedure, we first focus on a single stock as an example. Figure \ref{fig:spg} reports the stock price for Simon Property (SPG) throughout the detection period and its corresponding e-process initiated on January 1, 2007. Observing from the e-process figure, it is evident that both the e-mixture method and the e-GREE method effectively reject the null hypothesis at thresholds of 2, 5, 10, and 20 before the financial crisis ends. Notably, the e-GREE method generally takes fewer trading days compared to the e-mixture method to achieve this rejection across various threshold levels. Also, the null hypothesis is rejected using e-GREE method prior to another significant decline in the stock price during February 2009 to June 2009, thus preventing potential larger losses and underscoring the effectiveness of e-process methods.


Compared to e-batch and other p-variable based methods stated in Section \ref{sec:multi}, e-process based methods exhibit a unique advantage in sequential settings, particularly in financial applications where actual losses accumulate sequentially over time. In such scenarios, the e-process permits the early termination without a specified sampling period, potentially preventing further losses at an earlier stage. 
Table \ref{tab:eprocess} displays the number of trading days required to reject the null hypothesis at various threshold levels for the selected 20 stocks from 10 different sectors and the two stocks in Real Estate. 
The table shows that stocks in sectors significantly impacted by the 2007--2008 subprime crisis, such as Financials and Consumer Discretionary, Energy 
could generally be detected using e-process based methods. In particular, the representative companies in Real Estate are rejected the earliest; see the last rows of Table \ref{tab:eprocess}. 
 In contrast, for stocks in sectors less affected by the subprime crisis, such as Technology, Health Care, and Consumer Staples, 
we are unable to reject the null hypothesis. This is intuitive, given that their prices and returns remain relatively stable or even increase during the   financial crisis.

\begin{table}[t]
\caption{The number of trading days taken to detect evidence against $H(\hat{\mu},\hat{\sigma})$ using the e-GREE method and the e-mixture method for different stocks from January 1, 2007 to December 31, 2010; “–” means no detection is observed till December 31, 2010.}
\renewcommand{\arraystretch}{0.7}
\centering
\resizebox{\linewidth}{!}
{
\tabcolsep=0.2cm
\begin{tabular}{cccccccccccc}
\toprule
& & & \multicolumn{4}{c}{E-GREE} & & \multicolumn{4}{c}{E-mixture} \\
\cmidrule{4-7}\cmidrule{9-12} & Threshold &       & 2     & 5     & 10    & 20    &       & 2     & 5     & 10    & 20 \\
\toprule
\multirow{2}[2]{*}{Financials} & Bank of America  &       & 378	& 385 & 385	&  393&	& 	393	& 394& 	395	& 403 \\
& Morgan Stanley &       & 429  & 	439	  & 445  & 447	&   &  447 & 447	  & 447  & 	447\\
\midrule
\multirow{2}[2]{*}{Utilities} & The Southern  &       & -     & -     & -     & -     &       & -     & -     & -     & - \\
          & Duke Energy &       & -     & -     & -     & -     &       & -     & -     & -     & - \\
    \midrule
    Communication & Verizon Comms.   &       & -     & -     & -     & -     &       & -     & -     & -     & - \\
    Services      & AT\&T   &       & -     & -     & -     & -     &       & -     & -     & -     & - \\
    \midrule
Consumer    & Walmart  &    & -&	-	&-&	-	&	& -	&-	&-	&-   \\
  Staples   & PepsiCo &       & - &	- & - & -	& & 	-	& -& 	-& 	-  \\
\midrule
Consumer    & Ford Motor  &    & 476&	491	&498&	565	&	& 546	&594	&594	&594   \\
  Discretionary   & Las Vegas Sands &       & 442 &	445 & 447 & 450	& & 	451	& 454& 	457& 	457  \\
\midrule   
    
    \multirow{2}[2]{*}{Energy} 
& Texas Pacific Land &  &  158	& 244& 	261& 	269	& & 	242	& 261	& 261& 	263    \\
& Pioneer 
          & & 496&622&-&-&&-&-& -&- \\
\midrule

\multirow{2}[2]{*}{Material} 
& Southern Copper &       & 476	&496&	537&-& & 539 &- &- & - \\
& Air Products 
          &       &476   & 516   &  537 & -     &       & -     & -     & -     & - \\
\midrule

\multirow{2}[2]{*}{Health Care} & Johnson \& Johnson & & -   & -   & -   & -  &  & -  & -     & -     & - \\
          & Pfizer &  & -   & -  & -   & - &  & -     & -     & -     & - \\             
\midrule
\multirow{2}[2]{*}{Technology} & Int. Business Machines  & & -   & -   & -   & -  &  & -  & -     & -     & - \\
          & Microsoft  &  & -   & -  & -   & - &  & -     & -     & -     & - \\   
\midrule
\multirow{2}[2]{*}{Industrials} & General Electric  &  & 537	& 546 &	578&  - &   & -	  & -     & -     & - \\
& United Parcel Service & & 476	& 524& 	542	& 632 &  &542  &604  & - & - \\   
    \midrule    \midrule
    \multirow{2}[2]{*}{Real Estate } & Simon Property  &       & 165	& 224& 	242	& 254& & 		223& 	239& 	250	& 253\\
          & Prologis &      &264	&271	&271&	296	&&	270	&271	&271	&275\\             
    \bottomrule
    \end{tabular}
    }
  \label{tab:eprocess}
\end{table}


\section{Discussion}
\label{sec:7}

This paper proposes an e-process based approach for testing mean and variance from non-stationary data. 
We consider four classes of non-parametric composite hypotheses with specified mean and variance bound along with additional constraints of distribution, such as symmetry, unimodality, or a combination thereof. 
For this purpose, our main technical results  give the best p-variables and e-variables in the  simple setting where one summary data point is observed. The explicit formulas are summarized in Table \ref{tab:1}. 
Using the obtained e-variables, we construct an e-process using either the e-mixture method or the e-GREE method. 
Simulation studies and empirical analysis are conducted to show the performance of the proposed methods in comparison with GRAPA of \cite{WR24}  and with  the exponential supermartingale methods of  \cite{HRMS20, HRMS21} and \cite{WR23}.

As mentioned in Section \ref{sec:single}, our constructions of p-values and e-values are potentially useful for multiple testing, which is not addressed in this paper. 
The literature on using e-values in multiple testing is growing recently.
For instance, e-values are used for false discovery control in knockoffs; 
see  \cite{RB23} for derandomization, \cite{ALM23} for Bayesian linear models, and \cite{GS23} for resolution-adaptive variable selection. Finally, the obtained e-variables may also be useful to build e-confidence regions (see \cite{VW23})  and e-posterior as (see \cite{G23}) for  $(\mu,\sigma^2)$, although we mainly consider a non-parametric setting.

\section{Proofs of all results}\label{sec:5}
We collect all proofs in the paper in this section. 
\begin{proof}[Proof of Lemma \ref{lem:1}]
Let $\mathcal P$  be any collection of p-variables for $H$. 
For $Q\in H$, using the fact that the elements of $\mathcal P$ are comonotonic, we have
 $$Q(\inf \{P\in \mathcal P\}> \alpha) = Q\left(\bigcap_{P\in\mathcal P} \{ P>\alpha\} \right) = \inf_{P\in \mathcal P} Q (  P>\alpha   ) \ge 1- \alpha. $$
This implies 
 $$Q(\inf \{P\in \mathcal P\}\le \alpha) \le \alpha. $$
Hence, the infimum of all p-variables for $H$ is still a p-variable, which is the smallest one.   
\end{proof}

For all theorems below, we will prove precision statements for the formulation of $\E^Q[X]=\mu$ instead of $\E^Q[X]\le \mu$, making these statements stronger. 
For the validity statements, it is easy to verify that those p-variables and e-variables are valid under both formulations. 
\begin{proof}[Proof of Theorem \ref{th:1}]
Since the problem is invariant under location shift and scaling, it suffices to consider the normalized case of $(\mu,\sigma)=(0,1)$. 

It is clear that $P$ is decreasing in $X$ and $E$ is increasing in $X$.

For $Q\in H(0,1)$,
Cantelli's inequality implies $Q(X>x) \le 1/(1+x^2)$ for $x>0$, which implies, for each $\alpha \in (0,1)$, $$Q(P\le \alpha) = Q(1+X_+^2 \ge 1/\alpha) = Q\left( X\ge \sqrt{(1-\alpha)/\alpha} \right) \le \frac{1}{1+ {1/\alpha-1}} = \alpha.
$$
The inequality above is an equality if $Q$ is chosen such that  
\begin{align}
\label{eq:Qalpha-1}
Q\left (X=  \sqrt{(1-\alpha)/\alpha} \right )=\alpha=1-Q\left(X =- \sqrt{\alpha/(1-\alpha)}\right),
\end{align} 
and we can easily verify that $\E^Q[X]=0$ and $\var^Q(X)=1$.
This implies that $\sup_{Q\in H({0,1})} Q (P\le \alpha) =\alpha$ for each $\alpha\in (0,1)$, and therefore $P=1/(1+X_+^2) $ is a precise p-variable for $H(0,1)$.

For $Q\in H(0,1)$, we have $\E^Q[X_+^2] \le \E^Q[X^2]\le 1$.
To show that $E$ is precise, let $Q$ be given by \eqref{eq:Qalpha-1}, which satisfies $\E^Q[X_+^2] = \alpha $.
By taking $\alpha\uparrow 1$ we know $\sup_{Q\in H({0,1})} \E^Q[E]=1$, and therefore $E=X_+^2 $ is a precise e-variable for $H(0,1)$.  
\end{proof}

\begin{proof}[Proof of Theorem \ref{th:2}]
We first show the statement on the e-variable. Set $(\mu,\sigma)=(0,1)$ as in the proof of Theorem \ref{th:1}. 
For   $Q\in H _{\rm S}(0,1)$, we have $2 \E^Q[X_+^2] = \E^Q[X^2]\le 1$, with equal sign holding if $\var^Q(X)=1$. Therefore, $E=2 X^2_+$ is a precise e-variable for $H _{\rm S}(0,1)$.  

Since $E=2X^2_+$ is an e-variable, by Markov's inequality, 
 $1/E=(2 X_+)^{-2}$ is a p-variable for $H_{\rm S}(0,1)$. 
 In Theorem \ref{th:1} we have seen that $P_0$ is   a p-variable for $H (0,1)$, and hence also a p-variable for  $H_{\rm S}(0,1)\subseteq H (0,1)$. 
 Using Lemma \ref{lem:1}, the minimum of $P_0$ and $(2 E_0)^{-1}$ is a p-variable for $H_{\rm S}(0,1)$.

Next, we show that $P$ is semi-precise.
 For $\alpha \in (0,1/2]$, let $Q$ be chosen such that  
$$ 
Q\left(X=   (2\alpha) ^{-1/2} \right )=\alpha= Q\left(X =  -  (2\alpha) ^{-1/2}  \right) \mbox{~and~} Q(X=0)=1-2\alpha.
$$ 
We can  verify that $\E^Q[X]=0$, $\var^Q(X)=1$, and $X$ is symmetrically distributed.
It follows that $Q(P \le \alpha) = Q( X=   (2\alpha) ^{-1/2})=\alpha$. This implies that $\sup_{Q\in H_{\rm S}({0,1})} Q (P\le \alpha) =\alpha$ for $\alpha \in (0,1/2]$. 
Therefore, $P $ is a semi-precise p-variable for $H_{\rm S}(0,1)$. 

Finally, we show that there do  not exist precise p-variables for $H_{\rm S}(0,1)$. 
Suppose that $P=f(X)$ is a precise p-variable, where $f$ is a decreasing function.
Note that $Q(X\le 0)\ge 1/2$ for all $Q\in H_{\rm S}(0,1)$.
It follows that  $Q(P \ge  f(0)) \ge 1/2$ and $Q(P < f(0)) \le 1/2$. 
If $f(0)>1/2$, then for $\alpha \in[1/2, f(0)]$, $Q(P \le  \alpha ) \le 1/2<\alpha$, implying that $P$ is not precise. 
If $f(0)\le 1/2$,
then, by taking $Q$ as the point-mass at $0$, we have $Q(P\le 1/2)=1$, implying that $P$ is not a p-variable. Either way we have a contradiction, and hence does not exist a precise p-variable.
\end{proof}

\begin{proof}[Proof of Theorem \ref{th:3}]
Set $(\mu,\sigma)=(0,1)$ as in the proof of Theorem \ref{th:1}.  
By Theorem 1 of \cite{BKV20},  
\begin{align}\label{eq:BKV}
\sup_{Q\in H_{\rm U}(0,1)} T^Q_X(1-\alpha )  =  \max\left\{ \sqrt{\frac{4}{9\alpha }-1}  ,   \sqrt{\frac{3-3\alpha}{1+3\alpha}}\right\}  \mbox{~~~for $\alpha \in (0,1)$}.
\end{align}
Note that $P$ is a decreasing function of $X$, and we denote  this by $P=f(X) $ where $$ f (x)= \max\left\{ \frac 4 9(1+x_+^2 )^{-1},  \frac{4}{3}(1+ x _+^2 )^{-1}-\frac{1}{3}  \right\}.$$
For $\alpha \in (0,1/6]$, we have  $$\sup_{Q\in H_{\rm U}(0,1)}  T^Q_X(1-\alpha )=  \sqrt{\frac{4}{9\alpha }-1},$$ and hence
\begin{align*}
\inf_{Q\in H_{\rm U}(0,1)}  T^Q_P(\alpha )& = f\left(\sup_{Q\in H_{\rm U}(0,1)}  T^Q_X(1-\alpha)\right)
 =  \frac 4 9\left(1+  {\frac{4}{9\alpha }-1}  \right)^{-1} 
 =  \alpha.
\end{align*}  
For $\alpha \in ( 1/6,1)$, it is  
$$\sup_{Q\in H_{\rm U}(0,1)} T^Q_X(1-\alpha ) =  
 \sqrt{\frac{3-3\alpha}{1+3\alpha}},
 $$
and hence 
\begin{align*}
\inf_{Q\in H_{\rm U}(0,1)}  T^Q_P(\alpha )& = f\left(\sup_{Q\in H_{\rm U}(0,1)}  T^Q_X(1-\alpha)\right)
=
   \frac{4 }{3}\left (1+   {\frac{3-3\alpha}{1+3\alpha}} \right)^{-1}-   \frac{1}{3}  
  =   \alpha.
\end{align*}   
Using Lemma \ref{lem:2}, 
we obtain that $P$ is a precise p-variable for $H_{\rm U} (0,1)  $. 

As  $E$ is an e-variable for $H(0,1)$, it is also an e-variable for $H_{\rm U}(0,1)$. To show that it is precise, fix any $a\in (0,1)$, and let $p>0$ and $b>0$ satisfy
$$
a^2 = \frac{3-3p}{3p+2-p^2}  \mbox{~~and~~} b=\frac{1+p}{1-p} a.
$$
Note that such $p$ exists for any $a\in (0,1)$ since the range of $({3-3p})/({3p+2-p^2})$ covers $(0,1)$.
Choose $Q$ such that the distribution of $X$ has a point-mass at $-a$  with probability $p$
and a uniform density on $[-a, b]$.
We can compute
$$
\E^Q [X] =-a p + \frac{b-a}{2} (1-p) = -ap+ap= 0,
$$
and 
\begin{align*}
\E^Q [X^2 ] &= a^2 p + \frac {a^2}{3} (1-p) +\frac {b^2}{3} (1-p) 
 =   \frac{ a^2 (3p+2-p^2) }{3(1-p)} =1.
\end{align*}
Therefore $Q\in H_{\rm U}(0,1)$. 
We also have $$\E^Q[E]= \E^Q[X_+^2] = 1- a^2 p - \frac {a^2}{3} (1-p)\ge 1-a^2.$$
Since $a\in (0,1)$ is arbitrary, we get  $\sup_{Q\in H_{\rm U}(0,1)} \E^Q[E]=1$, and hence $E$ is a precise e-variable.
\end{proof}

\begin{proof}[Proof of Lemma \ref{lem:3}]
For $\alpha \ge 1/2$, since $Q\in H_{\rm US}(0,1)$ is symmetric about $0$,
we have $T^Q_X(1-\alpha)\le 0$, with $T^Q_X(1-\alpha)=0$ if $Q$ is the point-mass at $0$.
We assume $\alpha < 1/2$ below.


Take any $Q\in H_{\rm US}(0,1)$, and we will find another distribution $R$ with smaller variance and the same $\alpha$-quantile (we omit ``left" because the quantile is unique for $Q$ and $R$).
Note that $Q$ has  a decreasing density on $(0,\infty)$ and possibly a point-mass at $0$.
Denote by $x_0=T^Q_X(1-\alpha)$ and $g $ the density function of $Q$ on $(0,\infty)$.
Consider a different distribution $R$ symmetric with respect to $0$ which has uniform density equal to $g(x_0)$ on $(0,b)$ for some $b>x_0$ and a point-mass at $0$,
such that $R([x_0,b ))=\alpha=Q([x_0,\infty))=R([x_0,\infty))$. Denote by $h$ the density function of $R$ on $(0,\infty)$, and note that $h(x)=0$ for $x>b$. Since $Q$ has a decreasing density $g$ on $(0,\infty)$,
 $g \ge h$ on $(0,x_0)$
 and $g\le h$ on $(x_0,b)$.
The above conditions imply
\begin{align}\label{eq:g-h}
\int_0^{x_0} x^2  g(x)\d x \ge \int_0^{x_0} x^2  h(x)\d x 
 \mbox{~~and~~} 
 \int_ {x_0}^\infty x^2  g(x)\d x \ge \int_{x_0}^\infty x^2  h(x)\d x ,\end{align}
  where   the second inequality is due to $R([x_0,\infty))=Q([x_0,\infty))$.
  Note that both inequalities in \eqref{eq:g-h} are equalities if and only if $g=h$, and equivalently, $Q=R$. 
  It follows that $\E^Q[X^2] \ge  \E^R[X^2]$, and hence $R\in H_{\rm US}(0,1)$.
Note that the condition $Q([x_0,\infty))=\alpha=R([x_0,\infty))$ guarantees   $T^Q_X(1-\alpha)= x_0=T^R_X(1-\alpha)$; that is $R$ has the same $\alpha$-quantile as $G$.

The above argument shows that it suffices for us to consider distributions $Q$ which can be represented by a mixture of point-mass at $0$ and a uniform distribution on $[-b,b]$. 
We   also assume that $Q$ has variance $1$; if the variance is less than $1$, then a rescaled distribution from $Q$ has variance $1$ and a larger $\alpha$-quantile.  
Let $p=Q((0,\infty))\in (0,1/2]$. 
We can compute $\E^Q[X^2] = 2 p b^2/3  = 1$, and hence $b=3^{1/2} (2p)^{-1/2}$.
This gives $$T^Q_X(1-\alpha) =b (1-\alpha/p) = \sqrt{\frac{3}{p}} \left(1-\frac{\alpha}{p}\right).$$
Maximizing the above term over $p\in (0,1/2]$ gives $p=3\alpha$ if $\alpha \le 1/6$ and $p=1/2$ if $\alpha \in (1/6,1/2]$, showing the desired supremum formula in the lemma.
\end{proof}

\begin{proof}[Proof of Theorem \ref{th:4}] Set $(\mu,\sigma)=(0,1)$ as in the proof of Theorem \ref{th:1}. 
By Theorem \ref{th:3}, $E=2E_0$ is an e-variable for $H_{\rm US}(0,1)$. 
It is precise because  $\E^Q[2X_+^2]=1$ for any $Q\in H_{\rm US}(0,1)$ with $\var^Q(X)=1$.

The fact that precise p-variables do not exist for $H_{\rm US}(0,1)$ follows from the same argument as in the proof of the corresponding statement in Theorem \ref{th:2}. 

It remains to show that $P$ is a semi-precise p-variable for $H_{\rm US}(0,1)$.
Write $P=f(X) $ where
$$ f (x)= \frac{2}{9x^2} \id_{[ 4/3 ,\infty)}(x_+^2) +  \frac{3-\sqrt{3  }x }{6}\id_{ (0, 4/3 )}(x_+^2) + \id_{(-\infty,0]}(x).$$ 
Using Lemma \ref{lem:3}, for $\alpha \in (0,1/6]$, we have  $$\sup_{Q\in H_{\rm US}(0,1)}T^Q_X(1-\alpha) =  \sqrt{ \frac{2}{9\alpha} } \ge \sqrt{ \frac 43} $$ and 
\begin{align*}
\inf_{Q\in H_{\rm US}(0,1)} T^Q_P(\alpha )& = f\left(\sup_{Q\in H_{\rm US}(0,1)}T^Q_X(1-\alpha)\right)
 =  \frac 2 9 \times \frac{9\alpha}{ 2}
 =  \alpha.
\end{align*}  
Similarly, for  $\alpha \in (1/6,1/2)$, we have  $$\sup_{Q\in H_{\rm US}(0,1)}T^Q_X(1-\alpha) =  \sqrt 3{(1-2\alpha)}  \in \left(0, \sqrt{   4/3}\right)$$ and 
\begin{align*}
\inf_{Q\in H_{\rm US}(0,1)} T^Q_P(\alpha )& = f\left(\sup_{Q\in H_{\rm US}(0,1)}T^Q_X(1-\alpha)\right)
 =  \frac{3-3(1-2\alpha) }{6} 
 =  \alpha.
\end{align*}  
Finally, for $\alpha \in [1/2,1)$, we have 
$\inf_{Q\in H_{\rm US}(0,1)} T^Q_P(\alpha )=1$ since $\p(X\le 0)\ge 1/2$.
Using Lemma \ref{lem:2},
the above three cases together imply that 
$P$ is a semi-precise p-variable for $H_{\rm US}(0,1)$.
\end{proof}

\begin{proof}[Proof of a statement in Section \ref{sec:42}]
Here we show that p-Simes and p-Fisher can be applied to conditionally valid p-values. 
Assume $\mathbb{P} (P_t \leq \alpha | \mathcal F_{t-1}) \leq \alpha$ for each $t=1,\dots,n$ and $\alpha \in(0,1)$ under $H_0$. This implies that there exists $\tilde P_t\le P_t$ such that $\mathbb{P} (\tilde P_t \leq \alpha |\mathcal F_{t-1} ) = \alpha$ for all $\alpha \in (0,1)$. Hence, $\tilde P_1,\dots,\tilde P_n$ are iid. 
Applying the combination methods to $\tilde  P_1,\dots,\tilde  P_n$  yields a valid Type-I error control. Since $\tilde P_t\le P_t$ for each $t$ and the two combination methods are monotone, we also have a valid Type-I error control when combining $P_1,\dots,P_n$.  
\end{proof}

\begin{proof}[Proof of Proposition \ref{prop:consistency}]
The assumption that data are iid implies that $E_1,E_2,\dots$ are iid. 
The ``only if" statement is trivial since $\E^Q[E_1]\le 1$ implies that $(M_t)_{t\ge 1}$ is an e-process for $Q$, and hence $Q(\sup_{t\in [n]}M_t\ge 1/\alpha) \le \alpha$ for all $n\in \N$. 
Next we show the ``if" statement. 
For this, we use Theorem 3 of \cite{WWZ22}, which states that, under the iid assumption, 
$$\frac{1}{t} \left (\log M_{T}(\boldsymbol \lambda^{\rm GREE} ) - \log M_{t}(\boldsymbol \lambda^{\rm GRO} ) \right)\xrightarrow{L^1(Q)}0 \mbox{~~~
as $t\to\infty$},$$
where $M_{t}(\boldsymbol \lambda^{\rm GREE} )$ is given by \eqref{eq:e-process} 
with each $\lambda_i$ computed form the e-GREE method,
and  $M_{t}(\boldsymbol \lambda^{\rm GRO} )$ is given by \eqref{eq:e-process} 
with each $\lambda_i $ given by its theoretically growth-rate optimal value $$\lambda^*=\argmax_{\lambda \in (0,1]} \E^Q[\log (1-\lambda +\lambda E_1)],$$
and this gives $$\frac{1}{t} \log M_{t}(\boldsymbol \lambda^{\rm GRO} ) =\max_{\lambda \in (0,1]}   \E^Q[\log (1-\lambda +\lambda E_1)]. $$ 
Therefore, we have 
$$
\frac 1 t \log M_t \xrightarrow{Q} \max_{\lambda \in (0,1]}   \E^Q[\log (1-\lambda +\lambda E_1)] \mbox{~~as $t\to \infty$}.
$$
It remains to verify 
$\max_{\lambda \in (0,1]}\E^Q[\log (1-\lambda +\lambda E_1)]>1$.
Note that  $\E [E_1]>1 $
implies $\E [E_1\wedge K]>1$ for some $K\ge 1$.
We denote by $Y=E_1\wedge K$.
Since $\E  [(Y-1)_+] - \E  [(Y-1)_-] = \E  [Y-1]>0$, there exists some $\epsilon\in (0,1)$ such that  
$$\frac{1}{1+\epsilon}\E [(Y - 1)_+ ] - \frac{1}{1-\epsilon}\E [(Y  - 1)_- ] >0.$$
Note that
$\log (1+x) \ge x/(1+\epsilon)$ for $x\in [0,\epsilon)$
and $\log (1+x) \ge x/(1-\epsilon)$ for $x\in (-\epsilon,0)$, that is,
$$
\log(1+x) \ge \frac{x_+}{1+\epsilon}
-\frac{x_-}{1-\epsilon} \mbox{~~~for $x\in (-\epsilon,\epsilon).$}
$$
Hence, for $\lambda\in (0,\epsilon/K)$, implying $\lambda (Y-1)\in (-\epsilon,\epsilon)$, we have  
\begin{align*}
\E [\log (1-\lambda +\lambda E_1)]&\ge \E [\log (1 + \lambda (Y-1) )] \\& 
\ge   \frac{1}{1+\epsilon}\E [ \lambda (Y-1)_+ ] 
- \frac{1}{1-\epsilon}\E [ \lambda (Y-1)_- ]  
>0,
\end{align*}  
thus showing the desired inequality. 
\end{proof}

\begin{proof}[Proof of Proposition \ref{prop:power}]
First, it is clear that $ \Pi^{Q_\sigma}\ge 0$ by choosing $\lambda =0$  in the supremum. 
Second, by Jensen's inequality, for $\sigma>1$, 
 \begin{align*} \E^{Q_0}[\log ( 1-\lambda + \lambda \sigma^2 E)] 
& \le    \log ( 1-\lambda + \lambda \sigma^2 \E^{Q_0}[E])
 \le  \log ( 1-\lambda + \lambda \sigma^2 ) =2  \log \sigma.
\end{align*}
We next show $ \Pi^{Q_\sigma} \ge 2  p   \log \sigma  -\log 2    $. Note that
 \begin{align*} \E^{Q_0}[\log ( 1-\lambda + \lambda \sigma^2 E)] 
& \ge (1-p)    \log ( 1-\lambda   ) +p     \log ( 1-\lambda + \lambda \sigma^2  ).
\end{align*}
Maximizing the above term over $\lambda \in [0,1]$, 
the maximizer is 
$\lambda ^* = (p\sigma^2 -1)/(\sigma^2-1)$.
The corresponding maximum value satisfies
\begin{align*}
  (1-p)    \log \frac{(1-p)\sigma^2}{\sigma^2 -1} +p     \log ( p \sigma^2  )
&  \ge (1-p) \log (1-p) + p \log p + p\log \sigma^2
\\&  \ge - \log 2 + p\log \sigma^2,
\end{align*}
where we used the fact that $x\log x + (1-x)\log (1-x)$ on $[0,1]$ is maximized at $x=1/2$.
This shows  $ \Pi^{Q_\sigma} \ge  2 p   \log \sigma  -\log 2    $, completing the proof of \eqref{eq:rw-r1-1}.

Finally, we prove the last statement  $ 0\le \Pi^{Q_\sigma} -\Pi^{Q_\delta} \le 2(\log \sigma -\log \delta )$ for $\sigma>\delta>1$.
For any $\lambda \in [0,1]$, let $\lambda' =\lambda \delta^2 /\sigma^2 \in [0,1]$. We have 
$$ \Pi^{Q_\sigma} 
\ge \log (1-\lambda'+\lambda' \sigma^2 E) \ge \log( 1-\lambda +\lambda \delta^2 E).
$$
Taking a supremum over $\lambda\in [0,1]$ yields $ \Pi^{Q_\sigma} \ge \Pi^{Q_\delta} $.
To show the other inequality, 
\begin{align*}
\Pi^{Q_\sigma} 
&\le  \sup_{\lambda \in [0,1]}\log\left (\frac{\sigma^2}{\delta^2}(1-\lambda)+\lambda \sigma^2 E\right) 
\\&= \log \frac{\sigma^2}{\delta^2} +  \sup_{\lambda \in [0,1]}\log (1-\lambda+\lambda \delta^2 E) = 2 \log \frac{\sigma}{\delta} +   
\Pi^{Q_\delta}.
\end{align*}
This gives $ \Pi^{Q_\sigma} -\Pi^{Q_\delta} \le 2(\log \sigma -\log \delta )$ and completes the proof. 
\end{proof}

 \subsection*{Acknowledgements}
We thank the Editor, an Associate Editor, and two anonymous referees for constructive comments. 
We also thank Aaditya Ramdas, Qiuqi Wang, and Ian Waudby-Smith for helpful discussions.   Wang was partly supported by the Natural Sciences and Engineering Research Council of Canada.

\end{document}